\documentstyle[amscd]{amsart}
\newtheorem{thm}{Theorem}[section]
\newtheorem{lem}[thm]{Lemma}
\newtheorem{cor}[thm]{Corollary}
\newtheorem{prop}[thm]{Proposition}

\newtheorem{th}{Theorem}

\theoremstyle{definition}

\def \CPb {\overline{\bold{CP}}^{\,2}}
\def \z {\zeta}
\def \R {\bold{R}}
\def \Z {\bold{Z}}
\def \Sig{\Sigma}

\def \CM {\cal M}

\def \ASD {anti-self-duality }
\def \asd {anti-self-dual }

\def \la {\langle}
\def \ra {\rangle}
\newtheorem{conj}{Conjecture}

\def \a {\alpha}
\def \b {\beta}
\def \g {\gamma}
\def \d {\delta}
\def \k {\kappa}
\def \lam {\lambda}
\def \G {\Gamma}
\def \l {\ell}
\def \o {\omega}
\def \s {\sigma}
\def \t {\tau}
\def \z {\zeta}
\def \ba {\bar{\alpha}}
\def \bd {\partial}
\def \x {\times}
\def \ve {\varepsilon}
\def \D {\bold{D}}
\def \hD {\hat{D}}
\def \K {\bold{K}}
\def \A {\bold{A}}
\def \Q {\bold{Q}}
\def \bk {\bar{\k}}
\def \e {\epsilon}

\def \bL {\bar{L}}

\begin{document}

\baselineskip.525cm

\title[Rational Blowdowns]{Rational Blowdowns of Smooth $4$-Manifolds}
\author[Ronald Fintushel]{Ronald Fintushel}
\address{Department of Mathematics, Michigan State University \newline
\hspace*{.375in}East Lansing, Michigan 48824}
\email{ronfint@@math.msu.edu}
\thanks{The first author was partially supported NSF Grant DMS9401032 and the
second
author by NSF Grant DMS9302526}
\author[Ronald J. Stern]{Ronald J. Stern}
\address{Department of Mathematics, University of California \newline
\hspace*{.375in}Irvine,  California 92717}
\email{rstern@@math.uci.edu}

\maketitle

\section{Introduction\label{Intro}}

The invariants of Donaldson and of Seiberg and Witten are powerful tools for
studying
smooth $4$-manifolds. A fundamental problem is to determine procedures which
relate smooth $4$-manifolds in such a fashion that their effect  on both the
Donaldson
and Seiberg-Witten invariants can be computed.  The purpose of this paper is to
initiate
this study by  introducing  a surgical procedure, called  rational blowdown,
and to
determine  how this procedure affects these two sets of invariants. The
technique of
rationally blowing down and its effect on the the Donaldson invariant
were first announced at the 1993 Georgia International Topology Conference and
represents the bulk of the  mathematics in this paper. We fell upon this
surgical
procedure while we were investigating the behavior of the  Donaldson invariant
in the
presence of embedded spheres and while investigating methods for producing
a topological logarithmic transform.  As it turns out, this rational blowdown
procedure
allows for the full computation of the Donaldson series (and Seiberg-Witten
invariants) of all elliptic surfaces with $p_g \ge 1$ with the only input being
the
Donaldson invariants of the Kummer surface; in particular this computation
shows that
the Donaldson series of elliptic surfaces is that conjectured by Kronheimer and
Mrowka
in \cite{KM}:
\begin{th} Let $E(n;p,q)$ be the simply connected elliptic surface with
$p_g=n-1$ and with multiple fibers of relatively prime orders $p,q\ge1$. Then
\[{\D}_{E(n;p,q)}=\exp(Q/2){\sinh^n(f)\over\sinh(f_p)\sinh(f_q)}.\]
\end{th}
\noindent This theorem gives another, more topological, proof of the
diffeomorphism classification of elliptic surfaces
(\cite{Bauer,MorganMrowka,MorganOGrady,Fried1}). This procedure also goes
further and
routinely computes the Donaldson series (and Seiberg-Witten invariants) for
many
$4$-manifolds, some of which are complex surfaces, and for most of the
currently known
examples which are not even homotopy equivalent to complex surfaces.

The ideas presented in this paper have led to rather easy proofs of the blowup
formulas
for the Donaldson invariants for arbitrary  smooth $4$-manifolds
\cite{FSblowup} and
alternate proofs and generalizations
\cite{FSstructure} of some of the results announced by Kronheimer and Mrowka
(\cite{KM},\cite{KMbigpaper}). While we chose to first write up these later
results,
another major delay in the appearance of this paper was the introduction of the
Seiberg-Witten invariants.

{ }From the beginning, Witten has conjectured  how the Seiberg-Witten
invariants and the
Donaldson invariants determine each other (cf. \cite{Witten}). Some progress in
proving this relationship has been announced by V. Pidstrigach and A. Tyurin.
Our
techniques verify Witten's conjecture for elliptic surfaces and for a large
class of
manifolds obtained from them by rational blowdowns. (See \S 8.)

Here is an outline of the paper:  In \S 2  we introduce the
concept of a rational blowdown and discuss relevant topological issues. Our
main analytical result,  Theorem~\ref{basic},  gives a universal formula which
relates
the Donaldson invariants of a manifold with those of its rational blowdown.
Three
examples of the effect of a rational blowdown are  given in \S 3 and these
examples are
used in subsequent sections to compute the universal quantities given in
Theorem~\ref{basic}.  In  \S 4 we give the fundamental definitions of the
Donaldson
series, and \S 5 presents our key analytical results. Here we shall take
advantage of
our later results and techniques (\cite{FSblowup},\cite{FSstructure}) to
streamline
our earlier arguments. In particular, we will utilize the  \lq\lq \ pullback
--- pushforward "
point of view introduced and  developed by Cliff Taubes in
\cite{Sxl,Reds,Circle,Holo}
(or, alternatively the thesis of Wieczorek \cite{W}) to prove our basic
universal formula
(Theorem~\ref{basic}).  Under the assumption of simple type, this universal
formula
takes on a particularly simple form (Theorem~\ref{BASIC}). Starting with the
computations of the Donaldson series for elliptic surfaces without multiple
fibers
given in \cite{KM},\cite{FSstructure} and \cite{L}, we  apply
Theorem~\ref{BASIC} and
some of the  examples presented in \S 3 to compute the Donaldson series of the
elliptic
surfaces with multiple fibers in \S 6. Under the assumption of simple type and
the additional assumption that the configuration of curves that is blown down
is
`taut',  Theorem~\ref{BASIC} yields a very simple formula relating the basic
classes of $X$ with those of its rational blowdowns (cf.
Theorem~\ref{tautcalc}). This,
as well as applications to the  computations of the Donaldson series of other
manifolds,
is  discussed in \S 7.  Theorem~\ref{BASIC} has a straightforward analogue
relating the
Seiberg-Witten invariants of $X$ and those of its rational blowdowns. We
conclude this
paper with a statement and proof of this relationship in \S 8.

\bigskip

\section{The Topology of Rational Blowdowns\label{topology}}

In this section we  define what is meant by a rational blowdown.  Let $C_p$
denote the
simply-connected smooth
$4$-manifold obtained by plumbing the $(p-1)$ disk bundles  over the $2$-sphere
according to the linear diagram

\centerline{\unitlength 1cm
\begin{picture}(5,2)
\put(.9,.7){$\bullet$}
\put(1,.8){\line(1,0){1.3}}
\put(2.2,.7){$\bullet$}
\put(2.3,.8){\line(1,0){.75}}
\put(3.3,.8){.}
\put(3.5,.8){.}
\put(3.7,.8){.}
\put(4,.8){\line(1,0){.75}}
\put(4.65,.7){$\bullet$}
\put(.35,1.1){$-(p+2)$}
\put(2.1,1.1){$-2$}
\put(4.55,1.1){$-2$}
\put(.45,.4){$u_{p-1}$}
\put(2.1,.4){$u_{p-2}$}
\put(4.55,.4){$u_1$}
\end{picture}}

\noindent Here, each node denotes a disk bundle over $S^2$  with Euler class
indicated
by the label; an interval  indicates that the endpoint disk bundles are
plumbed,  i.e
identified fiber to base over the upper hemisphere  of each $S^2$. Label the
homology
classes represented by the spheres in $C_p$ by $u_1,\dots,u_{p-1}$ so that the
self-intersections are
$u_{p-1}^2=-(p+2)$ and, for $j=1,\dots,p-2$, $u_j^2=-2$. Further,  orient the
spheres
so that $u_j\cdot u_{j+1}=+1$.  Then $C_p$ is a  $4$-manifold with negative
definite
intersection  form and with boundary the lens space $L(p^2,p-1)$.

\begin{lem}\label{ratball} The lens space $L(p^2,p-1)=\partial C_p$ bounds a
rational
ball
$B_p$ with \
$\pi_1(B_p)={\Z}_p$ and a surjective inclusion induced homomorphism \
$\pi_1(L(p^2,p-1)={{\Z}}_{p^2}\to \pi_1(B_p)$.
\end{lem}
\begin{pf} There are several constructions of $B_p$; we present three here.
The first construction is perhaps amenable to showing that if  the
configuration of
spheres $C_p$ are symplectically embedded in a symplectic
$4$-manifold $X$, then the rational blowdown $X_p$ is also symplectic (cf.
\cite{Gompf}). For this construction let
${\bold{F}}_{p-1}$, $p\ge 2$, be the simply connected ruled surface  whose
negative
section $s_-$ has square $-(p-1)$. Let $s_+$ be a positive section (with square
$(p-1)$) and $f$ a fiber.  Then the homology classes
$s_++f$ and $s_-$ are represented by embedded $2$-spheres which intersect each
other
once  and have intersection matrix
\[ \begin{pmatrix} p+1& 1\\ 1 & -(p-1) \end{pmatrix} \] It follows that the
regular
neighborhood of this pair of $2$-spheres has boundary
$L(p^2,p-1)$. Its complement in ${\bold{F}}_{p-1}$ is the rational ball $B_p$.

The second construction begins with the configuration of $(p-1)$ $2$-spheres

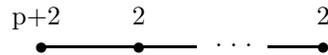
\centerline{\unitlength 1cm
\begin{picture}(5,2)
\put(.9,.7){$\bullet$}
\put(1,.8){\line(1,0){1.3}}
\put(2.2,.7){$\bullet$}
\put(2.3,.8){\line(1,0){.75}}
\put(3.3,.8){.}
\put(3.5,.8){.}
\put(3.7,.8){.}
\put(4,.8){\line(1,0){.75}}
\put(4.65,.7){$\bullet$}
\put(.6,1.1){p+2}
\put(2.2,1.1){2}
\put(4.65,1.1){2}
\end{picture}}

\noindent in $\#(p-1){\bold{CP}}^{\,2}$ where the spheres (from left to right)
represent
\[ 2h_1-h_2+\cdots-h_{p-1}, \ h_1+h_2, \ h_2+h_3, \dots , h_{p-2}+h_{p-1}\]
where $h_i$
is the hyperplane class in the $i$\,th copy of ${\bold{CP}}^{\,2}$. The
boundary of the
regular neighborhood of the configuration is $L(p^2,p-1)$ and the classes of
the
configuration span $H_2({\bold{CP}}^{\,2};{\Q})$. The complement is the
rational ball
$B_p$.

The third construction is due to Casson and Harer \cite{CH}. It utilizes the
fact that any
lens space is the double cover of $S^3$ branched over a 2-bridge knot. The
2-bridge knot
$K((1-p)/p^2)$ corresponding to $L(p^2,1-p)$ is slice, and $B_p$ is the double
cover of the
$4$-ball branched over the slice disk. \end{pf}

That all these constructions produce the same rational ball $B_p$ is an
exercise in Kirby calculus.
However, for the purposes of this paper, it is the third construction that is
the
most useful, since it allows us to quickly prove:

\begin{cor} Each diffeomorphism of $L(p^2,1-p)$ extends over the rational ball
$B_p$.
\label{lensdiff}\end{cor}
\begin{pf} It is a theorem of Bonahon \cite{Bonahon} that
$\pi_0(\text{Diff}(L(p^2,1-p))=\Z_2$, and is generated by the deck
transformation $\t$ of
the double branched cover of $K((1-p)/p^2)$. The extension of $\t$ to $B_p$ is
given by the
deck transformation of the double cover of $B^4$ branched over the slice
disk.\end{pf}

Suppose that $C_p$ embeds in a closed smooth $4$-manifold $X$. Then let $X_p$
be the
smooth $4$-manifold obtained by removing the interior of $C_p$ and replacing it
with
$B_p$. Corollary~\ref{lensdiff} implies that this construction is well-defined.
We call this
procedure a {\bf rational blowdown} and say that
$X_p$ is obtained by {\bf rationally blowing down} $X$. Note that
$b^+(X)=b^+(X_p)$ so that
rationally blowing down increases the signature while keeping $b^+$ fixed. An
algebro-geometric analogue of rationally blowing down is discussed in
\cite{KSB}.

With respect to the  basis $\{u_1,\dots,u_{p-1}\}$ for $H_2(C_p)$, the plumbing
matrix
for $C_p$ is given by the  symmetric $(p-1)\times(p-1)$ matrix
\[ P= \begin{pmatrix} -2 & 1& & & & &  \\ 1 & -2& 1& & &0 &  \\ 0 & 1&-2 &1 & &
&  \\
  & & & &\ddots & &  \\
 & 0& & & & -2& 1 \\ & & & & &1 & -(p+2)
\end{pmatrix}\]  with inverse given by $(P^{-1})_{i,j} =-j+{(ij)(p+1)\over
p^2}$ for
$j\le i$.

Let $Q:H_2(C_p,\partial C_p;{{\Z}})\x H_2(C_p;{\Z}) \to {\Z}$ be the (relative)
intersection  form of $C_p$ and  let $\{\g_1,\dots,\g_{p-1}\}$ be the basis of
$H_2(C_p,\partial C_p;\Z)$ dual  to the basis $\{u_1,\dots, u_{p-1}\}$ of
$H_2(C_p;{\Z})$
with respect to $Q$. I.e.
$\gamma_k\cdot u_\l=\delta_{k\l}$. Let $i_\ast:H_2(C_p;{{\Z}}) \to
H_2(C_p,\partial
C_p;{\Z})$  be the inclusion induced homomorphism. Then the intersection form
of
$H_2(C_p,\partial C_p;{\Q})$ is defined by
\[ \g_k\cdot\g_\l={1\over p^2}\g_k \cdot {\g}'_\l \]  where  ${\g}'_\l\in
H_2(C_p;{\Z})$ is chosen such that
$i_\ast({\g}'_\l)=p^2\g_\l$.   Since ${\g}'_\l=p^2P^{-1}({\g}_\l)$, the
intersection matrix
for $H_2(C_p,\partial C_p;{\Q})$ is $(\g_k\cdot\g_\l)=P^{-1}$.  Note also that
using
the sequence
\[\begin{CD} 0 \to H_2(C_p;{\Z}) @>P>> H_2(C_p,\partial C_p;{\Z}) @>\bd>>
H_1(L(p^2,1-p;{\Z})
\to 0 \end{CD}\]  we may identify $H_1(L(p^2,1-p;{\Z})$ with ${{\Z}}_{p^2}$ so
that
$\bd$  is given by $\bd(\g_j)=j$.

There is an alternative choice of dual bases for $H_2(C_p;{\Z})$ and
$H_2(C_p,\partial C_p;{\Z})$ that we shall find useful because of its symmetry.
Define
the basis $\{v_i\}$ of $H_2(C_p;{\Z})$ by
\[ v_i= u_{p-1}+\cdots + u_i, \hspace{.25in} u_j=v_j-v_{j+1} \] so
$v_i^2=-(p+2)$ for
each $i$, and if $i\ne j$ then $v_i\cdot v_j=-(p+1)$. The dual  basis
$\{\d_i\}$ of
$H_2(C_p,\bd C_p;{\Z})$ is given in terms of $\{\g_i\}$ by
\begin{eqnarray*} \d_i&=&\g_i-\g_{i-1}, \ \ i\ne 1\\
                  \d_1&=&\g_1  \end{eqnarray*}  Then
 \begin{eqnarray*} \d_i\cdot\d_j&=&{(p+1)\over p^2}, \ \ i\ne j\\
                  \d_i^2&=&-{(p^2-p-1)\over p^2}  \end{eqnarray*}  and
\[ \bd(\sum a_i\d_i) = \sum a_i. \]

Let the character variety of $SO(3)$ representations of
$\pi_1(L(p^2,1-p))$ mod conjugacy be denoted by $\chi_{SO(3)}(L(p^2,1-p))$, and
identify $\pi_1(L(p^2,1-p))$ with ${\Z}_{p^2}$ as above. Then we have an
identification
\[ \chi_{SO(3)}(L(p^2,1-p))\cong {{\Z}}_{p^2}/\{\pm1\}\cong
H_1(L(p^2,1-p);{{\Z}})/\{\pm1\}. \] Let $\eta$ be the generator of
$\chi_{SO(3)}(L(p^2,1-p))$ satisfying
\[ \eta(1) =\begin{pmatrix} \cos({2\pi i/p^2}) &\sin({2\pi i/p^2})&0\\
-\sin({2\pi i/p^2})& \cos({2\pi i/p^2})&0\\ 0&0&1  \end{pmatrix} \] Let $e\in
H_2(C_p,\bd C_p;{\Z})$; so $\bd e$ is some $n_e\in {\Z}_{p^2}$. Since
$b^+(C_p)=0$, $e$ defines an \asd connection $A_e$ on the complex line bundle
$L_e$ over
$C_p$ whose first chern class is the Poincar\'e dual of $e$. Throughout this
paper we
shall identify $H_2(C_p,\bd C_p;{\Z})\equiv H^2(C_p;{\Z})$; so we may write
$c_1(L_e)=e$.  Consider $C_p$ with a metric which gives a collar $L(p^2,1-p)\x
[
0,\infty)$. The connection $A_e$ has an asymptotic value as $t\to\infty$, and
this is a
flat connection on $L(p^2,1-p)$. Dividing out by gauge equivalence, we obtain
the
element
$\bd A_e=\eta^{n_e}\in\chi_{SO(3)}(L(p^2,1-p))$. For later use, we define
\[ \bd':H_2(C_p,\bd C_p;{{\Z}})\to\chi_{SO(3)}(L(p^2,p-1))
={{\Z}}_{p^2}/\{\pm1\}=\{0,1,\dots,[p/2]\} \] by $\bd'(e)=\bar{n}_e$, the
equivalence
class of $\bd e$.

\bigskip

\section{Examples of Rational Blowdowns\label{examples}}

In this section we present four examples of the effect of rational blowdowns.
These are
essential for our later computations.

\noindent {\bf Example 1.} Logarithmic transform as a rational blowdown

This first  example, whose discovery motivated our interest in this procedure,
shows
that  a logarithmic transform of order $p$ can be obtained by a sequence of
$(p-1)$
blowups (i.e. connect sum with $(p-1)$ copies of $\CPb$) and one rational
blowdown of
a natural embedding of the configuration $C_p$.
 First, some terminology. Recall that simply connected elliptic surfaces
without
multiple fibers are classified up to diffeomorphism by their holomorphic Euler
characteristic $n=e(X)/12=p_g(X)+1$. The  underlying smooth  $4$-manifold is
denoted
$E(n)$. The tubular neighborhood of a torus fiber  is a copy of $T^2\times
D^2=S^1\times(S^1\times D^2)$.  By a {\it log transform} on $E(n)$  we mean the
result
of removing this  $T^2\times D^2$ from $E(n)$ and regluing it by a
diffeomorphism
$$\varphi: T^2\times\partial D^2\to  T^2\times\partial D^2.$$ The {\it order}
of the
log transform is the absolute value of the degree of
$${\text {pr}}_{\partial D^2}\circ\varphi:{\text {pt}}\times \partial D^2\to
\partial
D^2.$$ Let $E(n)_{\varphi}$ denote the result of this operation on $E(n)$. Note
that
multiplicity 0 is a possibility.  It follows from  Moishezon
\cite{Moish} that if $\varphi$ and $\varphi'$ have the same order, there is a
diffeomorphism, fixing the boundary, from  $E(n)_{\varphi}$ to
$E(n)_{\varphi'}$. What
is needed here is the existence of a cusp neighborhood (cf. \cite{FScusp}). Let
$E(n;p)$  denote  any $E(n)_{\varphi}$ where the multiplicity of $\varphi$ is
$p$.

In $E(n;p)$ there is again a copy of the fiber $F$, but there is also a new
torus
fiber, the {\em multiple fiber}.  Denote its homology class by $f_p$; so in
$H_2(E(n;p);{\bold{Z}})$  we have $f=p\,f_p$. We can continue this process on
other torus
fibers; to insure that the  resulting manifold is simply connected we can take
at most
two  log-transforms with orders that are pairwise relatively prime.
 Let the orders be $p$ and $q$ and  denote the result by $E(n;p,q)$.  We
sometimes
write $E(n;p,q)$ in general, letting $p$ or $q$ equal $1$ if there are fewer
than $2$
multiple fibers. Of course one can take arbitrarily many log transforms (which
we
shall sometimes do) and we denote the result of taking $r$ log transforms of
orders
$p_1,\dots,p_r$ by $E(n;p_1,\dots,p_r)$.

The homology class $f$ of the fiber  of $E(n)$ can be represented  by an
immersed
sphere with one positive double point (a nodal fiber). Figure 1 represents a
handlebody
(Kirby calculus) picture for a cusp neighborhood $N$ which contains this nodal
fiber.
(See \cite{Kirby} for an explaination of such pictures and how to manipulate
them.)
Blow up this double point (i.e. take the proper  transform of $f$) so that the
class
$f-2e_1$  (where $e_1$ is the homology class of the exceptional divisor)  is
represented by an embedded sphere with square $-4$ (cf. Figure 2). This is just
the
configuration $C_2$.
 Now the exceptional divisor intersects this  sphere in two positive points.
Blow up
one of  these points, i.e. again take a proper transform.  One obtains the
homology
classes $u_2=f-2e_1-e_2$ and $u_1=e_1-e_2$ which form the configuration
$C_3$. Continuing in this fashion, $C_p$ naturally embeds in
$N\#_{p-1}{\CPb}\subset E(n)\#_{p-1}{\CPb}$ as in Figure 3.
Our first important example of a rational blown down is:

\begin{thm}\label{lgtr} The rational blowdown of the above configuration
$C_p\subset
E(n)\#(p-1)\CPb$ is diffeomorphic $E(n;p)$.
\end{thm}

\begin{pf} As proof, we offer a  sequence of Kirby calculus moves in Figures 4
through
8. In Figure 4 we add to Figure 3  the handle (with framing $-1$) which has
the property that when added to $\partial C_p$ one obtains $S^2\times S^1$ (so
that when a
further $3$ and $4$-handle are attached $B_p$ is obtained). Then we
 blow down the added handle, keeping track of the dual  2-handle (which is
labelled in
Figure 4  with 0-framing). In Figure 5 we blow down  this added handle with
framing
$-1$ and rearrange to obtain Figure 6. Now slide $e_1$ over the handle with
framing
$+1$ and rearrange to obtain Figure 7. Blow down the $-1$ curve in Figure 7; so
the $-2$ curve becomes a $-1$ curve. Continue this process $p-2$ times to
obtain Figure
8. If in this  final picture one replaces the handle with a dot on  it by a
1-handle,
there results the handlebody picture given by Gompf in \cite{nuc}  for $N_p$,
the order
$p$ log-transformed  cusp neighborhood. \end{pf}
\noindent For the case $p=2$, this  theorem was first observed by Gompf
\cite{Gompf}.

Here is a useful observation: To perform a log transform of order $pq$, first
perform
a log transform of order $p$ and then perform a log transform of order $q$ on
the
resulting multiple fiber $f_p$. This can also be obtained via a rational
blowdown
procedure.  Figure 9 is a  handlebody decomposition $N_p\#_{q-1}{\CPb}$ with an
easily
identified copy of $C_q$. The proof that  the result of blowing down $C_q$
results in
$E(n;pq)$ is to again follow through the steps of the proof of
Theorem~\ref{lgtr}.

\begin{prop}\label{ponq} Let $f_p$ be the multiple fiber in $E(n;p)$. Then
there is an
immersed (nodal) 2-sphere $S\subset E(n;p)$ representing the homology class of
$f_q$. Let $q$ be a positive integer relatively prime to $p$. If the process of
Theorem~\ref{lgtr} is applied to $S$, i.e. if $Y$ is the rational blowdown of
the
configuration $C_q$ in $E(n;p)\#(q-1)\CPb$ obtained from blowing up $S$, then
$Y\cong E(n;pq)$, the result of a multiplicity $pq$ log transform on
$E$.\end{prop}

\noindent {\bf Example 2.} In $E(2)$ there is an embedded sphere with
self-intersection $-4$ such that its blowdown is diffeomorphic to
$3{\bold{CP}}^2\#18\CPb$.

For this, any $-4$ curve suffices; however to verify that the rational blowdown
decomposes requires more Kirby calculus manipulations. The Milnor fiber
$M(2,3,5)$ for the Poincar\'e homology $3$-sphere $P=\Sigma(2,3,5)$ embeds in
$E(2)$ so
that $E(2)=M(2,3,5)\cup W$ for some $4$-manifold $W$ (cf. \cite{FScusp}). Now
$\partial M(2,3,5)=P$ also bounds another negative definite $4$-manifold $S$
which is the
trace of $-1$ surgery on the left handed trefoil. It is known that $S\cup W$ is
diffeomorphic to $3{\bold{CP}}^2\#11\CPb$.  Thus, to construct the example, it
suffices find a $-4$ curve in $M(2,3,5)$  whose rational blowdown produces
$S\#7\CPb$.
Recall that
$M(2,3,5)$ is just the $E_8$ plumbing manifold given in Figure 10. Slide the
handle
labeled $h$ over the handle labeled $k$ to obtain the $-4$ curve $h+k$ in
Figure 11.
Blow down this $-4$ curve to obtain Figure 12. Now slide the handle labeled
$h'$ over
the handle labeled $k'$ to obtain Figure 13. Now succesively blow down the $-1$
curves
to obtain Figure 14. Cancelling the $1-$ handle with the $2-$handle with
framing $-2$
 yields $S\#7\CPb$.

\noindent{\bf Example 3.} Given any smooth $4$-manifold $X$, there is an
embedding of
the configuration
$C_p\subset X\#(p-1)\CPb=Y$ with $u_i=e_{p-(i+1)}-e_{p-i}$ for $i=1,\dots,p-2$,
and
$u_{p-1}=-2e_1-e_2-\cdots-e_{p-1}$ such that the rational blowdown $Y_p$ of $Y$
 is
diffeomorphic to
$X\#H_p$ where $H_p$ is the homology $4$-sphere with $\pi_1={\Z}_p$ which is
the
double of the rational ball $B_p$.

In fact $C_p\subset \#(p-1){\CPb}=Y$, and, from the proof of
Lemma~\ref{ratball}, the
result  of blowing down this configuration is just the double of $B_p$.

Note that Example 3 points out that although a smooth $4$-manifold $Y$ may have
a
symplectic structure, it need  not be the case that  a rational blowdown $Y_p$
of $Y$
also have a symplectic structure. For in this example $X\#H_p$ will never have
a
symplectic structure since its $p-$fold cover can be written as a connected sum
of two
$4$-manifolds with positive $b_+$ so has vanishing Seiberg-Witten invariants
and hence,
by Taubes \cite{TSymplectic1}, is not symplectic.  Of course, in this example
the
configuration $C_p$ is not symplectically embedded. This  brings up the
possibility
that any smooth $4$-manifold can be rationally blown up to  a symplectic
$4$-manifold.

\bigskip

 \section{The Donaldson Series\label{def}}

In this section we outline the definition of the Donaldson invariant. We refer
the
reader to
\cite{Donpoly} and \cite{DK} for a more complete treatment. Given an oriented
simply
connected $4$-manifold with a generic Riemannian metric and an $SU(2)$ or
$SO(3)$ bundle
$P$ over $X$, the moduli space of gauge equivalence classes of \asd connections
on $P$
is a manifold ${\CM}_X(P)$ of dimension \[8\,c_2(P)-3\,(1+b_X^+)\] if $P$ is an
$SU(2)$
bundle, and
\[-2p_1(P)-3\,(1+b_X^+)\] if $P$ is an $SO(3)$ bundle. It will often be
convenient to
treat these two cases together by identifying ${\CM}_X(P)$ and
${\CM}_X(\text{ad}(P))$
for an $SU(2)$ bundle $P$. Over the product ${\CM}_X(P) \times X$ there is a
universal
$SO(3)$ bundle
${\bold{P}}$ which gives rise to a homomorphism $\mu:H_i(X;\R)\to
H^{4-i}({\CM}_X(P);\R)$
obtained by decomposing the class
$-{1\over 4}p_1({\bold{P}})\in H^4({\CM}_X \times X)$.

When either $w_2(P)\ne 0$ or when $w_2(P)= 0$, $d>\frac34(1+b_X^+)$, the
Uhlenbeck
compactification $\overline{{\CM}}_X(P)$ carries a fundamental class. In
practice, one is
able to get around this latter restriction by blowing up
$X$ and considering bundles over $X\#\CPb$ which are nontrivial when restricted
to the
exceptional divisor \cite{MMblowup}. In
\cite{FMbook} it is shown that for $\a\in H_2(X;{{\Z}})$ the classes
$\mu(\a)\in
H^2({{\CM}}_X(P))$ extend over $\overline{{\CM}}_X(P)$. When $b_X^+$ is odd,
$\dim
{\CM}_X(P)$ is even, say equal to $2d$. In fact, a class $c\in H_2(X;{{\Z}})$
and a
nonnegative integer
$d\equiv -c^2+\frac12(1+b^+)$ determine an $SO(3)$ bundle $P_{c,d}$ over $X$
with
$w_2(P_{c,d})\equiv c$ (mod $2$) and formal dimension $\dim {\CM}_X(P_{c,d}) =
2d$. For
$\ba=(\a_1,\dots,\a_d)\in H_2(X;{{\Z}})^d$, write
$\mu(\ba)=\mu(\a_1)\cup\cdots\cup
\mu(\a_d)$. Then one has
\[ \la
\mu(\ba),[\overline{{\CM}}_X(P_{c,d})]\ra=\int_{\overline{{\CM}}_X(P_{c,d})}\mu(\ba) \]
when $\mu(\ba)$ is viewed as a $2d$-form.

If $[1]\in H_0(X;{{\Z}})$ is the generator, then
$\nu=\mu([1])=-\frac14p_1(\b)\in
H^4({\CM}_X(P))$ where $\b$ is the basepoint fibration
$\tilde{{\CM}}_X(P)\to{\CM}_X(P)$ with
$\tilde{{\CM}}_X(P)$ the manifold of \asd connections on $P$ modulo based gauge
transformations, i.e. those that are the identity on the fiber over a fixed
basepoint.
The class $\nu$ extends over the Uhlenbeck compactification
$\overline{{\CM}}_X(P)$ if
$w_2(P)\ne0$, and in case $P$ is an $SU(2)$ bundle, the class will extend under
certain
dimension restrictions. Once again, these restrictions can be done away with
via the
tricks mentioned above \cite{MMblowup}.

Consider the graded algebra
\[{\A}(X)=\text{Sym}_*(H_0(X)\oplus H_2(X))\] where $H_i(X)$ has degree
$\frac12(4-i)$.
The Donaldson invariant $D_c=D_{X,c}$ is then an element of the dual algebra
${\A}^*(X)$, i.e. a linear function \[ D_c:{\A}(X)\to {\R}. \] This is a
homology
orientation-preserving diffeomorphism invariant for manifolds $X$ satisfying
$b_X^+\ge3$. Throughout this paper we assume $b_X^+\ge3$ and odd.

We let $x\in H_0(X)$ be the generator $[1]$ corresponding to the orientation.
In case
$a+2b=d>\frac34(1+b^+_X)$ and $\a\in H_2(X)$, \[
D_c(\a^ax^b)=\la\mu(\a)^a\nu^b,[{\overline{{\CM}}_X(P_{c,d})}]\ra\, . \]
We may extend
$\mu$ over ${\A}(X)$, and write for $z\in{\A}(X)$ of degree $d$,
$D_c(z)=\la\mu(z),[{\overline{{\CM}}_X(P_{c,d})}]\ra$. Since such moduli spaces
${\CM}_X(P_{c,d})$ exist only for $d\equiv -c^2+\frac12(1+b^+_X)$ (mod 4), the
Donaldson
invariant $D_c$ is defined only on elements of ${\A}(X)$ whose total degree is
congruent to
$-c^2+\frac12(1+b^+_X)$ (mod $4$). By definition, $D_c$ is $0$ on all elements
of other
degrees. When $P$ is an $SU(2)$ bundle one
simply writes $D$ or $D_X$.

If $Y$ is a simply connected $4$-manifold with boundary, one can similarly
construct
relative Donaldson invariants. A good reference for this is \cite{MMR}. When
the
boundary is a lens space, the theory simplifies considerably, and we get
relative Donaldson invariants
\[D_{Y,c}[\lam_i]:{\A}(Y)\to {\R}.\]

Following \cite{KM}, one considers the invariant \[
\hD_{X,c}:\text{Sym}_*(H_2(X))\to
{\R} \] defined by
$\hD_{X,c}(u)=D_{X,c}((1+\frac{x}{2})u)$. Whereas $D_{X,c}$ can be nonzero only
in
degrees congruent to $-c^2+\frac12(1+b^+)$ (mod $4$), $\hD_{X,c}$ can be
nonzero in degrees
congruent to $-c^2+\frac12(1+b^+)$ (mod 2). The {\em Donaldson series}
${\D}_c={\D}_{X,c}$ is defined by
\[{\D}_{X,c}(\a)=\hD_{X,c}(\exp(\a))=\sum_{d=0}^{\infty}{{\hD_{X,c}}(\a^d)\over
d!}\]
for all $\a\in H_2(X)$. This is a formal power series on $H_2(X)$.

A simply connected $4$-manifold $X$ is said to have {\em simple type} if the
relation
$D_{X,c}(x^2\,z)=4\,D_{X,c}(z)$ is satisfied by its Donaldson invariant for all
$z \in
{\A}(X)$ and for all $c\in H_2(X;{\Z})$. This important definition is due to
Kronheimer
and Mrowka \cite{KM} and was observed to hold for many $4$-manifolds
\cite{KMbigpaper,FSstructure}. In terms of
$\hD_{X,c}$, the simple type condition is that $\hD_{X,c}(zx)=2\hD_{X,c}(z)$
for all
$z\in {\A}(X)$ and for all $c\in H_2(X;{\Z})$. The assumption of simple type
assures
that for each $c$, the complete Donaldson invariant $D_{X,c}$ is determined by
the
Donaldson series ${\D}_{X,c}$.  It is still an open question whether all
$4$-manifolds are of simple type.

The structure theorem is:

\begin{thm}[Kronheimer and Mrowka
\cite{KMbigpaper,FSstructure}]\label{KMstruct}
Let $X$ be a simply connected 4-manifold of simple type. Then, there exist
finitely
many `basic' classes $\k_1$, \dots,
$\k_p\in H_2(X,{\Z})$ and nonzero rational numbers
$a_1$, \dots, $a_p$ such that \[{\D}_X\ =\ \exp(Q/2)\,\sum_{s=1}^p
a_se^{\k_s}\] as
analytic functions on $H_2(X)$. Each of the `basic classes' $\k_s$ is
characteristic,
i.e.
$\k_s\cdot x \equiv x\cdot x$ (mod $2$)for all $x\in H_2(X;{\Z})$.

Further, suppose $c\in H_2(X;{\Z})$. Then
\[ {\D}_{X,c}\ =\ \exp(Q/2)\,\sum_{s=1}^p(-1)^{{c^2+\k_s\cdot
c\over2}}a_se^{\k_s}\]
\end{thm}

\noindent Here  the homology class $\k_s$  acts on an arbitrary homology class
by
intersection, i.e. $\k_s(u)=\k_s\cdot u$. The basic classes $\k_s$ satisfy
certain inequalities analogous to the adjuction formula in a complex surface
\cite{KMbigpaper,FSstructure}. We shall need

\begin{thm}[\cite{FSstructure}]\label{FSadj}
 Let $X$ be a simply connected 4-manifold of simple type and let $\{\k_s\}$ be
the set
of basic classes as above. If $u\in H_2(X;{\Z})$ is represented by an immersed
$2$-sphere with $p\ge 1$ positive double points, then for each $s$
\begin{equation}
2p-2\
\ge u^2 + |\k_s\cdot u|. \label{adjintro} \end{equation}\end{thm}

\begin{thm}[\cite{FSstructure}]\label{FSadjspecial} Let $X$ be a simply
connected
4-manifold of simple type with basic classes
$\{\k_s\}$ as above. If the nontrivial class $u\in H_2(X;{\Z})$ is represented
by an
immersed
$2$-sphere with no positive double points, then let \[ \{\k_s|\,
s=1,\dots,2m\}\] be the
collection of basic classes which violate the inequality (\ref{adjintro}). Then
$\k_s\cdot u=\pm u^2$ for each such $\k_s$. Order these classes so that
$\k_s\cdot
u=-u^2\,(>0)$ for
$s=1,\dots,m$. Then
\[\sum_{s=1}^ma_se^{\k_s+u}-(-1)^{1+b_X^+\over2}\sum_{s=1}^ma_se^{-\k_s-u}=0.\]
\end{thm}

\bigskip

\section{The Basic Computational Theorem\label{gauge}}

Recall that for
$y\in H_2(X)$ and $F\in {\A}(X)$, interior product
\[ \iota_uF(v)= (\deg(v)+1) F(uv) \] defines a derivation which we denote by
$\bd_u$
and call `partial derivation'. Our basic theorem is:
\begin{thm}\label{basic} Let $X$ be a simply connected $4$-manifold of simple
type
containing the configuration $C_p$, and let $X_p$ be the result of rationally
blowing
down $C_p$. Then, restricted to $X^*=X_p \setminus B_p= X\setminus C_p$:
\[{\D}_{X_p}=\sum_{i=1}^{m(p)}\a_i(p)\bd^{n_i(p)}{\D}_{X,c_i(p)} \] where
$\a_i(p)\in
{\Q}$, $c_i(p)\in H_2(C_p;{\Z})$, $\bd^{n_i(p)}$ is an $n_i$th order partial
derivative
with  respect to classes in $H_2(C_p;{\Z})$, and these quantities depend only
on
$p$, not on $X$. \end{thm}

As motivation, and for use in the next section, we begin with a `by hand'
calculation.

\begin{lem}\label{C2} Let $X$ be a simply connected $4$-manifold containing an
embedded
2-sphere $\Sig$ of square $-4$ representing the homology class $\s$. Let $X_2$
be the
result of rationally blowing down $\Sig$. Then
\[ {{\D}_{X_2}|}_{X^*}={\D}_X-{\D}_{X,\s}. \]\end{lem}
\begin{pf} Here we work with $SU(2)$ connections over $X_2$ and $X$. The
conjugacy
classes of $SU(2)$ representations of $L(4,-1)$ are $\{\pm1,i\}$. Since a
multiple of
any class in $H_2(X_2;{\Z})$ lives in $H_2(X^*;{\Z})$, it suffices to evaluate
$D_{X_2}(z)$ for $z\in{\A}(X^*)$. The lemma is proved by a standard counting
argument
obtained by stretching the neck $\bd X^*\times {\R}$ in $X_2$. Doing this with
nonempty
moduli spaces leads to a sequence of
\asd connections (with respect to a sequence of generic metrics on $X_2$) which
limit
to \asd connections $A^*$ over $X^*$, and
$A_B$ over $B_2$ together perhaps with instantons on $X^*$ and $B_2$. Dimension
counting shows that $A^*$ is irreducible, $A_B$ is reducible (hence flat), and
that no
instantons occur. (The key fact is that each representation of $L(4,-1)$ has a
positive
dimensional isotropy group.) The flat $SU(2)$ connections on $B_2$ are
$\pm1$. Thus we have
\[ D_{X_2}(z)=\pm D_{X^*}[1](z)\pm D_{X^*}[-1](z). \] The invariants
$D_{X^*}[\pm1](z)$
are relative Donaldson invariants of $X^*$ with the given boundary values.

We first claim that $D_{X^*}[1](z)=\pm D_X(z)$. This is almost obvious by
applying an
argument like the one above. We need to know that there are no nontrivial
reducible
connections on the neighborhood $C_2$ of $\Sig$ with boundary value $1$ and in
a moduli
space of negative dimension. This follows simply from the fact that if $\lam$
is the
complex line bundle whose first chern class generates $H^2(C_2;{\Z})$, then the
moduli
space of \asd connections on $\lam^m+\lam^{-m}$ has dimension $4m-3$ (see
\cite{FSstructure}). To compute $D_{X^*}[-1](z)$, note that the Poincar\'e dual
of
$\s$ in $H^2(X;{\Z}_2)$ is the unique nonzero class whose restrictions to $X^*$
and
$C_2$ are both $0$. When passing to structure group $SO(3)$, the representation
$-1$
becomes trivial, and thus extends over $C_2$ as the trivial $SO(3)$ connection.
Now one
can see that $D_{X^*}[-1](z)=\pm D_{X,\s}(z)$.

Finally, we need to determine signs. A key point following from our discussion
is that
they are independent of $X$. Recall from Example 2 that there is a sphere
$\Sig$ of square $-4$ in the $K3$-surface $X$ which has a rational blowdown
$X_2$ with ${\D}_{X_2}=0$. Since ${\D}_{X,\s}=\exp(Q/2)={\D}_X$, our formula
must read
\[ D_{X_2}(z)=\pm(D_X(z)-D_{X,\s}(z)). \] To compute the overall sign, we must
compare
the way that signs are attached to
$A_0\#\Theta_{B_2}$,  and $A_0\#\Theta_{C_2}$ where $A_0$ is an \asd connection
on $X^*$
with boundary value $1$ and $\Theta_{B_2}$ and $\Theta_{C_2}$ are the trivial
connections on $B_2$ and $C_2$. This is done in a way similar to the proof of
\cite[Theorem 2.1]{FSblowup}, and the sign is easily seen to be `$+$'.
\end{pf}

We now proceed toward the proof of Theorem~\ref{basic}. The first step is to
understand
reducible connections over $C_p$. It will be convenient here to use the
symmetric dual
bases $\{v_i\}$ and $\{\d_i\}$ of  \S\ref{topology}. Using these coordinates,
we express
elements of  $H_2(C_p,\bd C_p;{\Z})$ as
\[ \b=\sum t_i\d_i = \la t_1,\dots,t_{p-1}\ra. \] Classes of the form $\la
t,\dots,t,s,\dots,s\ra $ will play a special role. We shall use the
abbreviation
\[ \la t,\dots,t,s,\dots,s \ra =\la t,s;b\ra \] if the number of $s$'s is $1\le
b\le
p-1$. If $e\in H_2(C_p,\bd C_p;{\Z})$, write
${\CM}_e$ for the $SO(3)$-moduli space of \asd connections on
$C_p$ which contains the reducible connection in the bundle $L_e\oplus {\R}$
where
$c_1(L_e)=e$, and which are asymptotically flat with boundary value $\bd' e\in
\chi_{SO(3)}(L(p^2,1-p))$. Note that $\bd\la t,t+1;b\ra =(p-1)t+b$.

\begin{lem}\label{dim} Let  $e=\la t,t+1;b\ra$ with $0\le t\le p$. Then \
$\dim{{\CM}}_e=2t-1$.
\end{lem}
\begin{pf} With respect to the basis $\{\d_i\}$, the intersection form of
$H_2(C_p,\bd)$  is
\begin{equation}\label{Q} Q= -{(p^2-p-1)\over p^2}\sum x_i^2 + 2\ {p+1\over
p^2}\sum_{i<j}x_ix_j\end{equation}   and
\begin{multline*} e^2=(b(t+1)^2+(p-b-1)t^2)(-{(p^2-p-1)\over p^2})\\
+2((p-b-1)bt(t+1)+
\binom{p-b-1}{2}t^2+\binom{b}{2}(t+1)^2){p+1\over p^2}\end{multline*}  Hence
\begin{equation}\label{esquare} e^2={1\over p^2}(b^2+b^2p-bp^2-2bt+t^2-pt^2).
\end{equation}  By hypothesis, $\bd e=(p-1)t+b\ne 0$. From \cite{Lawson} we
have
\[ {\rho\over 2}(\bd e)=-{1\over
p^2}(-2b^2-2b^2p-p^2+2bp^2+4bt-2p^2t-2t^2+2pt^2) \]
and by the index theorem \cite{APS}:
\[ \dim {\CM}_e= -2e^2-\frac32-\frac12(h+\rho)(\bd e)=-2e^2-2-{\rho\over 2}(\bd
e)=2t-1.\]
\end{pf}

\begin{lem}\label{bvlem1} Let $e=\la t,t+1;b\ra $ with $t\ge 0$ and
$(p-1)t+b\le p^2/2$.
Suppose also that
$e'=\la\a_1,\dots,\a_{p-1}\ra$ with $\sum\a_i=(p-1)t+b+rp^2$, $r\ne 0,-1$. Then
\
$\dim{\CM}_{e'}>\dim{\CM}_e$.\end{lem}
\begin{pf} Using \eqref{Q}, it follows from symmetry that for fixed $s=\sum
x_i$, the
minimum absolute value of
$Q\la x_1,\dots,x_{p-1}\ra$ occurs at $\mu(s)=\la s/(p+1),\dots,s/(p+1)\ra$,
and
\[ \mu(s)^2= -{(p^2-p-1)\over p^2}(p-1){s^2\over (p-1)^2} + 2\
{p+1\over p^2}\binom{p-1}{2} {s^2\over (p-1)^2} ={s^2\over p^2-p^3}.\] On the
other
hand by \eqref{esquare}, $e^2={1\over p^2}(b^2+b^2p-bp^2-2bt+t^2-pt^2)$.
Set $s= (p-1)t+b+rp^2$. Then
\[ \mu(s)^2-e^2=-{1\over p-1}(b+b^2+2br+p^2r^2+2rt(p-1)-bp).\] Since $1\le b\le
p-1$,
we have $bp\le p^2-p\le p^2r^2$. So
\[  \mu(s)^2-e^2\le -{1\over p-1}(b+b^2+2br+2rt(p-1))\] and if we assume $r\ge
1$, \
$\mu(s)^2< e^2 \text{(}<0\text{)}$. By the index theorem,
\[\dim {\CM}_{e'}= -2{e'}^2-\frac32-\frac12(h+\rho)(\bd e')\ge -2\mu(s)^2-
   \frac32-\frac12(h+\rho)(\bd e)\ge \dim {\CM}_e \]
since $(h+\rho)(\bd e')=(h+\rho)(\bd e)$. Notice that we have not yet used the
hypothesis that
$(p-1)t+b\le p^2/2$.

If $r<-1$, set $\bar{e}=\la t',t'+1;c\ra$ with $t',c$ chosen such that
\[ (p-1)t'+c=p^2-((p-1)t+b)\ge p^2/2.\]  By  Lemma \ref{dim},
$\dim{\CM}_{\bar{e}}\ge\dim{\CM}_e$ with equality only if $t'=t$.  Note that
$\dim{\CM}_{-e'}=\dim{\CM}_{e'}$, and $-\sum\a_i=(p-1)t'+c-(r+1)p^2$. Since
$-(r+1)\ge 1$, the case we have already handled shows that
$\dim{\CM}_{-e'}\ge\dim{\CM}_{\bar{e}}$.
\end{pf}

\begin{lem}\label{bvlem2} Let $e=\la t,t+1;b\ra $ with $t\ge 0$. Suppose that
$e'=\la\a_1,\dots,\a_{p-1}\ra \ne e$ but $\sum\a_i=(p-1)t+b$. Then
$\dim{\CM}_{e'}>\dim{\CM}_e$ unless $e'$ is a permutation of $e$.\end{lem}
\begin{pf} It suffices to show that ${e'}^2<e^2$. Write $e'=e+\nu$ where
\[\nu=\la n_1,\dots,n_{p-b-1},n_{p-b},\dots,n_{p-1}\ra.\] Since the sum of the
coordinates of $e$ and $e'$ is the same, $\sum n_i=0$. Let
\[ N_L=\sum_{i=1}^{p-b-1}n_i\hspace{.5in} N_R=\sum_{i=p-b}^{p-1}n_i.\]
\begin{eqnarray*} {e'}^2&=&e^2+2(N_L((p-2)t+b)+N_R((p-2)t+b-1))({p+1\over
p^2})\\
 && \hspace{2in} - 2(N_Lt+N_R(t+1))({p^2-p-1\over p^2}) +\nu^2\\
&=&e^2-2N_R+\nu^2
\end{eqnarray*} since $N_L+N_R=0$. Hence
$\frac12(\dim{\CM}_{e'}-\dim{\CM}_e)=e^2-{e'}^2=-\nu^2+2N_R$. However, if $y$
is the
result of adding $+1$ to $x_{i_0}$ and $-1$ to $x_{i_1}$ in
$x=\la x_1,\dots,x_{p-1}\ra$, then $y^2-x^2= 2(x_{i_1}-x_{i_0}-1)$. Starting
with
$x=\la0,\dots,0\ra$ and making these $\pm1$ moves with constant sign in each
coordinate
until reaching $\nu$, we see that the minimum change in the square is $-2$.
This is
achieved only if each coordinate operated on is originally $0$. Thus, if $N_+$
is the
sum of the positive coordinates $n_i$, we have  $-\nu^2\ge 2N_+$. Equality
occurs only
if each $n_i$ is $\pm1$ or $0$. In this case there are $N_+$ such $-1$'s. If
$|N_R|<N_+$
then $-\nu^2+2N_R\ge 2(N_+-|N_R|)>0$. If $|N_R|=N_+$ then each $-1$ occurs in a
coordinate $n_i$, $i= p-b,\dots,p-1$, and so $e'$ is a permutation of $e$. If
$-\nu^2>2N_+$ then since $|N_R|\le N_+$, we have $-\nu^2+2N_R>0$.
\end{pf}

\begin{prop}\label{bv} Let $e=\la t,t+1;b\ra $ with $t\ge 0$ and $(p-1)t+b\le
p^2/2$.
If $e'=\la\a_1,\dots,\a_{p-1}\ra$ with $e'\equiv e$ (mod $2$) and
$\dim{\CM}_{e'}\le\dim{\CM}_e$, then $\bd e' \le \bd e$ as elements of
${\Z}_{p^2}$.
\end{prop}
\begin{pf} Let $\bar{e}=\la s,s+1;c\ra$ with $s\ge 0$, be the unique class of
this form
with $0\le\bd\bar{e}\le p^2/2$ satisfying $\bd e'=\bd\bar{e}$. Lemmas
\ref{bvlem1} and \ref{bvlem2} imply that unless $-p^2/2\le\sum a_i<0$, we have
$\dim{\CM}_{\bar{e}}\le \dim{\CM}_{e'}$; so $s\le t$. This holds in any case,
since we
can always work with $-e'$. If $s=t$ then
$\bar{e}=e$ since no class $\la t,t+1;b'\ra$ with $b'\ne b$ is congruent to
$e\pmod2$.
This means that $\bd e' \le \bd e$.
\end{pf}

\begin{cor}\label{bddim} Let $e=\la t,t+1;b\ra $ with $t\ge 0$ and $(p-1)t+b\le
p^2/2$.
Suppose that $e'=\la\a_1,\dots,\a_{p-1}\ra$ with $\bd' e'=\bd'
e\in\chi_{SO(3)}(L(p^2,1-p))$ and
$e'\equiv e$ (mod $2$). Then $\dim{\CM}_{e'}=\dim{\CM}_e+4k$, $k\ge
0$.\end{cor}
\begin{pf} As above, $\dim{\CM}_e\le \dim{\CM}_{e'}$. But  $e'\equiv e$ (mod
$2$) implies
that ${e'}^2=e^2$ (mod $4$); so the corollary follows from the index
theorem.\end{pf}

We need one more simple fact. Let $\iota :(C_p,\emptyset)\to (C_p,\bd)$ be the
inclusion.
\begin{lem}\label{getc2} Let $e\in H_2(C_p,\bd;{\Z})$, and suppose that
$\bd e\equiv 0$ (mod $2$) in case $p$ is even. Then there is a \ $c\in
H_2(C_p;{\Z})$ such
that $\iota_*(c)\equiv e$ (mod $2$).\end{lem}
\begin{pf} This follows directly from the exact sequence
\[ 0\to H_2(C_p;{\Z})\to H_2(C_p,\bd;{\Z})\to {\Z}_{p^2}\to 0\, .\]
\end{pf}

We now proceed toward the proof of Theorem \ref{basic}. We shall work always
with
structure group $SO(3)$ and identify $SU(2)$ connections with $SO(3)$
connections on
$w_2=0$ bundles. We wish to calculate $D_{X_p}(z)$ for $z\in{\A}(X^*)$. If we
blow up
$X^*$ and evaluate $D_{X_p\#\CPb,e}(ze)=D_{X_p}(z)$ where $e$ is the
exceptional class
\cite{MMblowup}, we can work under the assumption that there are no flat
connections on
the complement of $B_p$ with the same $w_2$ as our given bundle. Keeping this
in mind,
we may simplify notation without loss by making the same assumption for our
given
situation, $X_p=X^*\cup B_p$. Consider a sequence of generic metrics on $X_p$
which
stretch a collar on $L(p^2,1-p)=\bd B_p$ to infinite length, giving the
disjoint union
of $X^*$ and $B_p$ with cylindrical ends as the limit. A sequence of
\asd connections $\{ A_n\}$ with respect to these metrics, each of which also
lies in the
divisor $V_z$ corresponding to $z\in{\A}(X^*)$, must limit to $A_{X^*} \amalg
A_{B_p}$.
These are \asd connections over $X^*$ and $B_p$, and a counting argument shows
that
$A_{X^*}\in V_z$ and $A_{B_p}$ is reducible. (Our above assumption is helpful
here.)
Since the only reducible connections on the rational ball $B_p$ are flat, we
get
\begin{equation}\label{bareqn} D_{X_p}(z)=\sum_{n=0}^{[p/2]}\pm
D_{X^*}[\eta^{np}](z).
\end{equation}  The notation $D_{X^*}[\eta^{np}]$ stands for the relative
Donaldson
invariant on
$X^*$ constructed from the moduli space of \asd connections over $X^*$ (with a
cylindrical end) which decay exponentially to a flat connection whose gauge
equivalence
class corresponds to the conjugacy class of the representation $\eta^{np}$.

We need to calculate the summands of \eqref{bareqn}. We begin with $n=0$, i.e.
$D_{X^*}[1](z)$. Consider $D_X(z)$. To calculate this, we use a neck-stretching
argument
as above. We see that on $C_p$ we must get a reducible \asd connection
corresponding to
chern class $e$ with $\dim{\CM}_e<0$ and
$e\equiv 0\pmod2$. This last condition means that $e$ cannot have the form
$\la 0,1;b\ra$ (recall $1\le b\le p-1$); so by Lemma \ref{dim}, $e\ne \la
t,t+1;b\ra$,
$t\ge 0$. Now Proposition \ref{bv} implies that $e=0$. Thus
\[ D_X(z)=\pm D_{X^*}[1](z) \]
and the sign is independent of $X$.

To calculate the other terms, we must utilize techniques of Taubes
\cite{Sxl,Reds,Circle,Holo} or Wieczorek \cite{W} as in
\cite[\S4]{FSstructure}. We
shall quickly review the methods involved and refer the reader to
\cite{FSstructure} and
the references given there for more details. Our plan is to evaluate all the
$D_{X^*}[\eta^m](z)$ inductively. (In case $p$ is even, we only need to
calculate this
for $m$ even.) We do this by computing
$D_{X,c_m}(z\,w_m)$ where $c_m\in H_2(X;{\Z})$ is supported in
$C_p$, $m=(p-1)t+b$, and $w_m\in \text{Sym}_t(H_2(C_p;{\Z}))$ depending only on
$m$ and
$p$. First we obtain $c_m$.  Let
$e_m\in H^2(C_p;{\Z})$ be the  Poincar\'e dual of $\la t,t+1;b\ra$.  By Lemma
\ref{getc2} we can find $c_m\in H_2(C_p;{\Z})\subset H_2(X;{\Z})$ such that
$\iota_*(c_m)\equiv \la t,t+1;b\ra\pmod2$. Thus the Poincar\'e dual of $c_m$ in
$H^2(X;{\Z})$ restricts to $C_p$ congruent to $e_m\pmod2$ and restricts
trivially to
$X^*$. A dimension counting argument shows that in the formalism of Taubes
\cite{Sxl},
$D_{X,c_m}(z\,w_m)$ is the sum of terms of the form
\begin{equation}\label{terms}
\int_{\tilde{{\CM}}_{X^*}[\eta^j]\x_j\tilde{{\CM}}_{C_p,\e,\l}}
    \tau\wedge\tilde{\mu}(z)\wedge\tilde{\mu}(w_m).  \end{equation} In this
formula,
$\tilde{{\CM}}_{C_p,\e,\l}$ is the based moduli space of exponentially decaying
asymptotically flat \asd connections on the $SO(3)$ bundle
$E_{\e,\l}$ which is obtained from the reducible bundle
$L_\e\oplus{\R}$ by grafting in $\l$ instanton bundles. (The euler class of
$L_\e$ is $\e$, $\bd\e = j$, $\e\equiv e_m\pmod2$, and
$\dim{\CM}_\e+8\l\le 2t-1$.) The notation `$\x_j$' in the formula denotes the
fiber
product with respect to the $SO(3)$-equivariant boundary value maps
\[ \bd_{C_p,\e,\l}:\tilde{{\CM}}_{C_p,\e,\l}\to G[j], \hspace{.25in}
  \bd_{X^*}[j]:\tilde{{\CM}}_{X^*}[\eta^j]\to G[j] \] where $G[j]\subset SO(3)$
is the
conjugacy class $\eta^j$ of representations of
$\pi_1(L(p^2,1-p))$ to $SO(3)$. If $j\ne 0,p^2/2$ then $G[j]$ is a 2-sphere,
$G[0]=\{ I\}$, and, in case $p$ is even, $G[p^2/2]\cong{\bold{RP}}^2$.  Also,
$\tau$
denotes a 3-form which integrates to 1 over the fibers of the basepoint
fibration
$\b_{X^*,j}$ i.e.
$\tilde{{\CM}}_{X^*}[\eta^j]\to {\CM}_{X^*}[\eta^j]$. The form
$\tilde{\mu}(w_m)$ is
supported near the orbit of the reducible connection corresponding to $\e$. (If
$\l>0$, this reducible connection lies in the Uhlenbeck compactification of
${\CM}_{C_p,\e,\l}$.) The principal $SO(3)$ bundle $\b_{X^*,j}$ has a reduction
to a
bundle with structure group $S^1$. As in \cite[\S4]{FSstructure}, we let
$\ve\in H^2({\CM}_{X^*}[\eta^j])$ denote the euler class of this $S^1$ bundle.

The upshot of Taubes' work cited above is that there is a form
$\tilde{\mu}(w_m)$ representing a class $\mu_{SO(3)}(w_m)$ in the
$SO(3)$-equivariant
cohomology of an enlargement of $\tilde{{\CM}}_{C_p,\e,\l}$. The lift
$\tilde{\mu}(z)$ defines an element of the equivariant cohomology
$H^{2d}_{SO(3)}(\tilde{{\CM}}_{X^*}[\eta^j])$. Furthermore, Taubes has shown
that the
push-forward $(\bd_{C_p,\e,\l})_*$ is well-defined, and
\[ \int_{\tilde{{\CM}}_{X^*}[\eta^j]\x_j\tilde{{\CM}}_{C_p,\e,\l}}
    \tau\wedge\tilde{\mu}(z)\wedge\tilde{\mu}(w_m)
   = \int_{\tilde{{\CM}}_{X^*}[\eta^j]}\tau\wedge\tilde{\mu}(z)\wedge
(\bd_{X^*}[j])^*(\bd_{C_p,\e,\l})_*(\tilde{\mu}(w_m)) \] where
$(\bd_{X^*}[j])^*$
denotes pullback in equivariant cohomology.

For $j=0$, $\bd_{C_p,\e,\l}:\tilde{{\CM}}_{C_p,\e,\l}\to \{1\}$, has fiber
dimension
equal to $\dim\tilde{{\CM}}_{C_p,\e,\l}= 4k+8\l$ for some
$k\ge 0$. The cohomology class of $(\bd_{C_p,\e,\l})_*(\tilde{\mu}(w_m))$ lies
in
$H^{2t-4k-8\l}_{SO(3)}(\{1\};{\R})=H^{2t-4k-8\l}(BSO(3);{\R})$ which is a
polynomial
algebra on the 4-dimensional class $\wp$, which pulls back over
$\tilde{{\CM}}_{X^*}[\eta^j]$ as $p_1(\b_{X^*,j})$. For $j\ne 0,p^2/2$, let
$j=t_j(p-1)+b_j$ where $1\le b_j\le p-1$. Then
$\bd_{C_p,\e,\l}:\tilde{{\CM}}_{C_p,\e,\l}\to G[j]$ has fiber dimension
$2t_j+2+8\l+4k-2$ for some $k\ge 0$; so the cohomology class of
$(\bd_{C_p,\e,\l})_*(\tilde{\mu}(w_m))$ lies in
$H^{2(t-t_j)-8\l-4k}_{SO(3)}(S^2;{\R})=
  H^{2(t-t_j)-8\l-4k}({\bold{CP}}^\infty;{\R})$. Let $v$ be the 2-dimensional
generator of
$H^*({\bold{CP}}^\infty;{\R})$. The pullback
$(\bd_{X^*}[j])^*(v)=\ve$. Using the fact that $\ve^2=p_1(\b_{X^*,j})$, and
arguing as
in \cite[Prop.4.5,4.6]{FSstructure} we get
\begin{equation}\label{expand}
 D_{X,c_m}(z\,w_m)=\sum_{t_j\equiv t\,(2)}\sum_q r_{m,j,q}D_{X^*}[\eta^j](zx^q)
   +\sum\begin{Sb} t_j\not\equiv t\,(2)\\j\ne 0\end{Sb}
       \sum_q r'_{m,j,q}D_{X^*}[\eta^j](zx^q\ve).
\end{equation} The notation $D_{X^*}[\eta^j](zx^q\ve)$ is not standard, but its
meaning
is clear. It follows from Proposition \ref{bv} that the $\eta^j$ in
\eqref{expand} have
$j\le m$; so this bounds $j$ in both terms. We emphasize that in order to
obtain
$r_{m,j,q}$ or $r_{m,j,q}'\ne 0$ we must have an $\e\in H^2(C_p;{\Z})$
satisfying
$\bd'\e = j$, $\e\equiv e_m\pmod2$, and $\dim{\CM}_\e+8q\le 2t-1$.

Assume inductively that:
\begin{itemize}
\item[a)]For each $j< m$ ($j\equiv 0\pmod2$ if $p$ is even)  there are classes
$w_{j,i}\in\text{Sym}_*(H_2(C_p;{\Z}))$ and rational numbers $a_{j,i}$
satisfying
\begin{equation}\label{induct}
    D_{X^*}[\eta^j](z)=\sum_{i=1}^j a_{j,i}D_{X,c_i}(zw_{j,i})
\end{equation}
\item[b)] For each $j$ with $t_j< t-1$ (and $j\equiv 0\pmod2$ if $p$ is even)
there are
classes  $w'_{j,i}\in\text{Sym}_*(H_2(C_p;{\Z}))$ and rational numbers
$a'_{j,i}$
satisfying
\begin{equation}\label{inductepsilon}
  D_{X^*}[\eta^j](z\ve)=\sum_{i=1}^j a_{j,i}D_{X,c_i}(zw'_{j,i})
\end{equation}
\end{itemize} for all $z\in{\A}(X^*)$, and the coefficients $a_{j,i},a_{j,i}'$
are
independent of $z$ and $X$.

Recall that we are writing $m=(t-1)p+b$ with $1\le b\le p-1$, and let $e_m$ be
the
Poincar\'e dual of $\la t,t+1;b\ra=(t+1)\g_{p-1}-\g_{p-1-b}$. Also, we suppose
that $m$
is even if $p$ is even. We set
\[ w_m=(u_{p-1}-(t-1)u_{p-1-b})\cdot (u_{p-1})^{t-1}\in{\A}(C_p). \] We wish to
calculate $D_{X,c_m}(z\,w_m)$ using
\eqref{expand}. For $j=m$ in this formula, we need to compute
$(\bd_{C_p,e_m,0})_*(\tilde{\mu}(w_m))\in H^0_{SO(3)}(G[m];{\R})={\R}$ since
$t_m=t$.
In fact,
\begin{multline*} (\bd_{C_p,e_m,0})_*(\tilde{\mu}(w_m))=r_{m,m,0}\\ =
-\frac12\la
u_{p-1}-(t-1)u_{p-1-b},e_m\ra\,(-\frac12\la u_{p-1},e_m\ra)^{t-1} =
(-\frac12)^t(2t)(t+1)^{t-1}\ne 0\end{multline*} (cf. \cite[p.187]{DK}).  In
\eqref{expand}, $ r_{m,m,0}D_{X^*}[\eta^m](z)$ is the only term which involves
the
boundary value
$\eta^m$. If $j$ is the boundary value of an $\e$ with  $\e\equiv e_m\pmod2$,
and
$\dim{\CM}_\e+8q\le 2t-1$, and if $t_j=t-1$, then by Corollary~\ref{bddim} and
Lemma~\ref{bvlem2}, $\e$ must be a permutation of $\la t-1,t;p-1-b\ra=
t\g_{p-1}-\g_b$.
In fact $\e\equiv e_m\pmod2$ implies that $\e = \la
t,t-1;b\ra=(t-1)\g_{p-1}+\g_{p-1-b}$. So
$j=(t-1)(p-1)+(p-1-b)$.  Hence
$\la u_{p-1}-(t-1)u_{p-1-b},\e\ra =0$. Thus, no such $j$ occurs in the second
sum of
the expansion \eqref{expand} for $D_{X,c_m}(z\,w_m)$. (I.e. for such
$j$, necessarily $q=0$ and $r'_{m,j,q}=0$.) Finally, if $p$ is even, then we
are
assuming that $m$ is also even. If $r_{m,i,q}$ or $r_{m,i,q}'\ne 0$ then as
above there
is an
$\e$ with $\bd\e = i$ and $\e\equiv e_m\pmod2$; so for
\[ \bd_2:H_2(C_p,\bd;{\Z}_2)\to H_1(L(p^2,1-p);{\Z}_2)={\Z}_2 \]
$j\equiv\bd_2(\e)\equiv\bd_2 e_m\equiv m\pmod2$. Accordingly, all the other
terms in
\eqref{expand} are given inductively by \eqref{induct} and
\eqref{inductepsilon}, and
the powers of $x$ can be removed using the hypothesis that $X$ has simple type.
Since the
coefficient of
$D_{X^*}[\eta^m](z)$ is nonzero, we may solve for it, completing the induction
step for
\eqref{induct}.

For \eqref{inductepsilon}, we show how to compute $D_{X^*}[\eta^{m'}](z\ve)$
for
$m'=(t-1)(p-1)+(p-1-b)$ as required. Thus after completing the inductive step
for each
$t(p-1)+c$, $1\le c\le p-1$, we will have completed the calculation of
$D_{X^*}[\eta^j](z\ve)$ for all $j=(t-1)p+a$, $1\le a\le p-1$. So to calculate
$D_{X^*}[\eta^{m'}](z\ve)$ and thus complete the induction, we calculate
$D_{X,c_m}(z\,w'_{m'})$ where
$w'_{m'}=(u_{p-1}+(t+1)u_{p-1-b})\cdot(u_{p-1}+(t-1)u_{p-1-b})\cdot
(u_{p-1})^{t-2}$.
Using \eqref{expand}
\begin{equation}\label{another} D_{X,c_m}(z\,w'_{m'})=\sum_{t_j\equiv
t\,(2)}\sum_q
s_{{m'},j,q}D_{X^*}[\eta^j](zx^q)
   +\sum\begin{Sb} t_j\not\equiv t\,(2)\\j\ne 0\end{Sb}
       \sum_q s'_{{m'},j,q}D_{X^*}[\eta^j](zx^q\ve).\end{equation} Computing as
above, we
see that $s_{m',m',0}=0$. What we need to see is that
$s'_{m',m',0}\ne 0$. By the argument of the above paragraph, $m'$ is the only
possible
boundary value not covered by the induction step.  Let
$\e=(t-1)\g_{p-1}+\g_{p-1-b}$. This is the only euler class that can give
boundary
value $m'$ in \eqref{another}. Then
\[ (\bd_{C_p,\e,0})_*(\tilde{\mu}(w'_m))\in H^2_{SO(3)}(S^2;{\R})\cong
H^2_{SO(3)}(G[m'];{\R})={\R} \] and  $(\bd_{C_p,\e,0})_*(\tilde{\mu}(w'_m))=
(-\frac12)^{t-1}2t(2t-2)(t-1)^{t-2}v$ which pulls back over
$\tilde{{\CM}}_{X^*}[\eta^j]$
as
$(-\frac12)^{t-3}t(t-1)(t-1)^{t-2}\ve$. This means that we can solve
\eqref{another} for
$D_{X^*}[\eta^{m'}](z\ve)$, completing the induction and the proof of Theorem
\ref{basic}.

The argument above shows that all of the relative invariants
$D_{X^*}[\eta^{np}]$ can be
expressed in terms of absolute invariants of $X$. Since we are assuming that
$X$ has
simple type, it follows that each of the relative invariants satisfies the
formula
\[ D_{X^*}[\eta^{np}](z\,x^2) = 4\,D_{X^*}[\eta^{np}](z). \] Hence it follows
from
\eqref{bareqn} that:

\begin{cor}\label{st} Let $X_p$ be the result of rationally blowing down
$C_p\subset X$.
If $X$ has simple type, then so does $X_p$. \ \ \qed \end{cor}

Now we shall make stronger use of the hypothesis that $X$ has simple type. By
\cite{KM,FSstructure} we can write
\begin{eqnarray*} {\D}_X&=&\exp(Q_X/2)\sum_{s=1}^na_se^{\k_s} \\
{\D}_{X,c}&=&\exp(Q_X/2)\sum_{s=1}^n(-1)^{\frac12(c^2+c\cdot\k_s)}a_se^{\k_s}
\end{eqnarray*}
for nonzero rational numbers $a_s$ and basic classes $\k_1,\dots,\k_n\in
H_2(X;{\Z})$. Here $Q_X$ is the intersection form of $X$. Now
\[ \bd_u(\exp(Q_X/2)e^{\k})=\exp(Q_X/2)(\tilde{u}+\k\cdot u)e^{\k} \]
where $\tilde{u}: H_2(X)\to\R$ is $\tilde{u}(\a)=u\cdot\a$ and
$\bd_v\tilde{u}=v\cdot u$.
Apply Theorem \ref{basic}: since all derivatives are taken with respect to
classes $u\in
H_2(C_p;\Z)$, after all derivatives are taken, the remaining $\tilde{u}$'s
restricted to
$X^*$ vanish. Hence,
\begin{equation}\label{barseries0}
{{\D}_{X_p}|_{X^*}}=\exp(Q_{X^*}/2)\sum_{s=1}^na_sb_se^{\k_s}|_{X^*}=
\exp(Q_{X^*}/2)\sum_{s=1}^na_sb_se^{\k'_s}
\end{equation}
where $\k'_s=\k_s|_{X^*}=\text{PD}(i^*(\text{PD}(\k_s)))\in
H_2(X^*,\bd;{\Z})$, where $\text{PD}$ denotes Poincar\'e duality, $i$ is the
inclusion $X^*\subset X$, and $b_s$ depends only on the intersection numbers of
$\k_s$
with the generators $u_i$ of $H_2(C_p;\Z)$.

\begin{lem}\label{theyextend} If $b_s\ne0$ in \eqref{barseries0} then
\[ \bd\k'_s\in p{\Z}_{p^2}\subset H_1(L(p^2,1-p);{\Z})={\Z}_{p^2}. \] \end{lem}
\begin{pf} Corollary~\ref{st} implies that $X_p$ has simple type. We thus have
\begin{equation}\label{barseries}
{\D}_{X_p}=\exp(Q_{X_p}/2)\sum_{r=1}^mc_re^{\lam_r}
\end{equation} where the basic classes of $X_p$ are $\lam_1,\dots,\lam_m$.
Restrict
${\D}_{X_p}$ to
$X^*$ and compare the restrictions of $\exp(Q_{X_p}/2)^{-1}{\D}_{X_p}$ in
\eqref{barseries0} and \eqref{barseries}. Since for distinct $\a\in
H_2(X^*,\bd;{\Z})$ the
functions $e^\a:H_2(X^*)\to{\R}$ are linearly independent, it follows that if
$b_s\ne 0$,
then $\k'_s=\lam_i|_{X^*}$ for some $i$. Thus $\k'_s$ extends over $B_p$, and
hence $\bd
\k'_s\in p{\Z}_{p^2}$. \end{pf} As a result, we have the following restatement
of
Theorem \ref{basic}.

\begin{thm}\label{BASIC} Suppose that $X$ has simple type and
\[{\D}_X=\exp(Q_X/2)\sum_{s=1}^na_se^{\k_s}.\]  Let $C_p\subset X$ and let
$X_p$ be its
rational blowdown. Let $\{\k_t|t=1,\dots,m\}$ be the basic classes of $X$ which
satisfy
$\bd\k'_t\in p{\Z}_{p^2}$, and for each $t$, let $\bk_t$ be the unique
extension of $\k'_t$.
Then
\[{\D}_{X_p}=\exp(Q_{X_p}/2)\sum_{t=1}^ma_tb_te^{\bk_t}\]  where the $b_t$
depend only
on the intersection numbers  $u_i\cdot\k_t$,
$i=1,\dots,p-1$.
\ \ \qed \end{thm}
\bigskip

\section{The Donaldson Invariant of Elliptic Surfaces\label{ellipticcompute}}

In this section we shall compute the result on the Donaldson series of
performing log
transforms.  The Donaldson invariants of the elliptic surfaces
$E(n), n\ge2$ without multiple fibers have been known for some time. There is a
complete calculation in \cite{FSstructure}, for example. For $n\ge2$:
\[ {\D}_{E(n)} = \exp(Q/2)\sinh^{n-2}(f) \] where $f$ is the class of a fiber.
In this
notation, the $K3$ surface is $E(2)$. As in Theorem~\ref{lgtr}, let
$X=E(2)\#(p-1)\CPb$, and let $X_p$ be the rational blowdown of $C_p\subset X$,
so that
$X_p\cong E(2;p)$.  Since ${\D}_{E(2)}=\exp(Q/2)$, the blowup formula
\cite{FSblowup}
yields
\begin{equation}\label{blowup} {\D}_X={1\over
2^{p-1}}\exp(Q/2)\sum_J\exp(\sum_{i=1}^{p-1}\e_{J,i}e_i)
\end{equation} where the outer sum is taken over all
$J=(\e_{J,1},\dots,\e_{J,p-1})\in\{\pm1\}^{p-1}$. The basic classes of $X$ are
$\{\k_J=\sum\e_{J,i}e_i\}$, and applying Theorem~\ref{BASIC} we get
\begin{equation}\label{DXbar}
{\D}_{X_p}={1\over2^{p-1}}\exp(Q_{X_p}/2)\sum_Jb_Je^{\bk_J}
\end{equation}
where $\bk_J\in H_2(X_p;{\Z})$ is the unique extension of
${\k_J|}_{X^*}$. Recall that the spheres of the configuration $C_p$ represent
homology
classes
$u_i=e_{p-(i+1)}-e_{p-i}$ for $1\le i\le p-2$, and
$u_{p-1}=f-2e_1-e_2-\cdots-e_{p-1}$.
In $X_p$ we have the multiple fiber $f_p=f/p$.

\begin{prop}\label{J} $\bk_J= |J|\cdot f_p$ where
$|J|=\sum_{i=1}^{p-1}\e_{J,i}$.\end{prop}
\begin{pf} First we find a class $\z\in H_2(C_p;{\Q})$ so that $(\k_J+\z)\cdot
u_i=0$
for each $i$. This means that $\k_J+\z\in H_2(X^*;{\Q})$, and as dual forms:
$H_2(X^*;{\Z})\to{\Z}$, ${\k_J|}_{X^*}=\k_J+\z$. To find $\z$ we need to solve
the
linear system
\[ (\k_J+\sum x_iu_i)\cdot u_j = 0,\ \ \ j=1,\dots,p-1. \] We begin by
rewriting
these equations. Let $\{\o_i\}$ be a standard basis for
${\Q}^{p-1}$, and let $A$ be the $(p-1)\x(p-1)$ matrix whose $i$th row vector
is
\begin{eqnarray*} A_i&=& \o_{p-(i+1)}-\o_{p-i}, \ \ i=1,\dots,p-2\\
                 A_{p-1}&=&-2\o_1-\o_2-\dots-\o_{p-1}. \end{eqnarray*} We have
$u_i=A^t(\o_i)\cdot{{\bold{e}}}$ and $u_{p-1}=f+A^t(\o_{p-1})\cdot{\bold{e}}$ ,
where
${\bold{e}}=(e_1,\dots,e_{p-1})$. Our linear system is equivalent to
\[ P{{\bold{x}}}=A\pmb{\e}_J\] where ${{\bold{x}}}=(x_1,\dots,x_{p-1})$ and
$\pmb{\e}_J=(\e_{J,1},\dots,\e_{J,p-1})$. (The matrix $P$ is the plumbing
matrix for
$C_p$.) Hence ${\bold{x}}=P^{-1}A\pmb{\e}_J$.

We claim that $P(A^t)^{-1}=-A$. This can be checked on the basis
\[ \{\o_2-\o_1,\dots,\o_{p-1}-\o_{p-2},\o_{p-1} \}\] using
\begin{eqnarray*} A(\o_i)&=&-\o_{p-1}-\o_{p-(i+1)}+\o_{p-i},\ 2\le i\le p-1 \ \
(\o_0=0),\\
A(\o_1)&=&-2\o_{p-1}+\o_{p-2},\\ P(\o_i)&=&\o_{i+1}-2\o_i+\o_{i-1},\ i\ne
p-1,\\
P(\o_{p-1})&=&-(p+2)\o_{p-1}+\o_{p-2}.
\end{eqnarray*}

It follows that $A^tP^{-1}A=-I$. Thus
\[\k_J+\z=\k_J+\sum
x_iu_i=({\bold{\e}}_J+A^t{{\bold{x}}})\cdot{\bold{e}}+x_{p-1}f
=(\pmb{\e}_J-\pmb{\e}_J)\cdot{\bold{e}}+x_{p-1}f=x_{p-1}f.\] To compute
$x_{p-1}$ note
that
\[
A\pmb{\e}_J=(\e_{J,p-2}-\e_{J,p-1},\e_{J,p-3}-\e_{J,p-2},\dots,\e_{J,1}-\e_{J,2},
-2\e_{J,1}-\e_{J,2}-\cdots-\e_{J,p-1}) \] so that if $(P^{-1})_{p-1}$ denotes
the
bottom row of $P^{-1}$:
\[x_{p-1}=(P^{-1})_{p-1}(A\pmb{\e}_J)=-{1\over
p^2}(1,2,\dots,p-1)\cdot(A{\bold{\e}}_J)
={1\over p}\sum\e_{J,i}={1\over p}|J|.\] Thus ${\k_J|}_{X^*}=\k_J+\z={1\over
p}|J|f$ as
forms: $H_2(X^*;{\Z})\to{\Z}$. The homology class $\k_J+\z$ is in fact an
integral
class $\bar{\k}_J=|J|f_p\in H_2(X_p;{\Z})$ which is the unique extension of
${\k_J|}_{X^*}$ \end{pf}

In an arbitrary smooth $4$-manifold $X$, define a {\em nodal fiber} to be an
immersed
2-sphere $S$ with one singularity, a positive double point, such that the
regular
neighborhood of $S$ is diffeomorphic to the regular neighborhood of a nodal
fiber in an
elliptic surface. (There need not be any associated ambient fibration of $X$.)
Given
such a nodal fiber $S$, one can perform a `log transform' of multiplicity $p$
by
blowing up to get $C_p\subset X\#(p-1)\CPb$ with
$u_{p-1}=S-2e_1-e_2-\cdots-e_{p-1}$,
and then blowing down $C_p$. We denote the result of this process by $X_p$.

Throughout, we use the following notation. If $X$ has simple type, and
$${\D}_X=\exp(Q/2)\sum a_se^{\k_s},$$ then we write ${\K}_X=\sum a_se^{\k_s}$.

\begin{prop}\label{formallog} Let $S$ be a nodal fiber which satisfies
$S\cdot\lam_j=0$
for each basic class $\lam_j$ of $X$. Then
\[{\D}_{X_p}=\begin{cases}
\exp(Q_{X_p}/2){\K}_X\cdot(b_{p,0}+\sum\limits_{i=1}^{p-1\over2}b_{p,2i}(e^{2iS/p}+e^{-2iS/p})),\
&p\ \text{odd}\\
\exp(Q_{X_p}/2){\K}_X\cdot(\sum\limits_{i=1}^{p\over2}b_{p,2i-1}(e^{(2i-1)S/p}+e^{-(2i-1)S/p})),
&p\ \text{even}\end{cases}\] where the coefficients $b_{p,j}$ depend only on
$p$, not on
$X$.
\end{prop}

\begin{pf} The Donaldson series of $X\#(p-1)\CPb$ is
\[
{1\over2^{p-1}}{\D}_X\cdot\exp(Q_{(p-1)\CPb}/2)\sum_J\exp(\sum_{i=1}^{p-1}\e_{J,i}e_i).
\]
Theorem~\ref{basic} states that ${\D}_{X_p}$ is obtained from this by applying
a
differential operator which by hypothesis evaluates trivially on ${\D}_X$. The
proposition now follows from \eqref{DXbar} and Propostion~\ref{J} by the
Leibniz rule.
(That the coefficients of $e^{mp}$ and $e^{-mp}$ are equal follows from the
fact that
${\D}_{E(2;p)}$ is an even function.) \end{pf}

\begin{prop}\label{log2} The Donaldson series of the simply connected elliptic
surface
$E(n;2)$ with $p_g=n-1$ ($>0$) and one multiple fiber of multiplicity $2$ is
\[ {\D}_{E(n;2)}=\exp(Q/2){\sinh^{n-1}(f)\over\sinh(f_2)}. \] \end{prop}
\begin{pf} According to Theorem~\ref{lgtr}, we obtain $E(n;2)$ from
$E(n)\#\CPb$ by
blowing down the sphere of square $-4$ representing $f-2e$. We have
${\D}_{E(n)\#\CPb}=\exp(Q/2)\sinh^{n-2}(f)\cosh(e)$. Lemma~\ref{C2} gives
\[ {{\D}_{E(n;2)}|}_{X^*} = ({\D}_{E(n)\#\CPb}-{\D}_{E(n)\#\CPb,f-2e})|_{X^*}=
       2\exp(Q/2)\sinh^{n-2}(f)\cosh(e)|_{X^*} \] (cf.\cite{KMbigpaper},
\cite[Thm.5.13]{FSstructure}). By Proposition~\ref{J}
\[
{\D}_{E(n;2)}=2\exp(Q/2)\sinh^{n-2}(f)\cosh(f_2)=\exp(Q/2){\sinh^{n-1}(f)\over\sinh(f_2)}.\]
\end{pf}

Proposition~\ref{formallog} now implies:
\begin{cor} If $S$ is a nodal fiber in $X$ orthogonal to all basic classes and
$X_2$ is
the multiplicity $2$ log transform of $X$ formed from $S$, then
\[{\D}_{X_2}=\exp(Q_{X_2}/2){\K}_X\cdot(e^{S/2}+e^{-S/2}). \ \ \qed\]
\end{cor}

\begin{lem}\label{sum} Let $X$ contain a nodal fiber $S$ orthogonal to all
basic
classes. Then the sum of the coefficients $b_{p,j}$ in the expression for
${\D}_{X_p}$ in
Proposition~\ref{formallog} is equal to $p$.\end{lem}
\begin{pf} In Example 3 we showed that there is a configuration
$C'_p\subset X\#(p-1)\CPb=Y$ where $u'_i=e_{p-(i+1)}-e_{p-i}$ for
$i=1,\dots,p-2$, and
$u'_{p-1}=-2e_1-e_2-\cdots-e_{p-1}$ such that the rational blowdown
$Y_p=X\#H_p$ where $H_p$ is a homology $4$-sphere with $\pi_1={\Z}_p$. It
follows
easily that ${\D}_{Y_p}=p\cdot {\D}_X$.

As above, we let $\k_J=\sum\e_{J,i}e_i$, $J\in\{\pm1\}^{p-1}$; so
\[ {\D}_Y={1\over2^{p-1}}{\D}_X\cdot\exp(Q_{(p-1)\CPb}/2)\sum_J e^{\k_J}. \]
All
partial derivatives of ${\D}_X$ with respect to classes in $H_2(C'_p)$ are
trivial;  so
\[ p{\D}_X={\D}_{{Y}_p}={\D}_X\cdot\sum_Jb_Je^{\bar{\k}_J}. \]
The proof of Proposition~\ref{J} shows that each $\bar{\k}_J=0$; so
$\sum_Jb_J=p$.

We can also form the configuration $C_p\subset Y$ whose blowdown is the $p$-log
transform of the nodal fiber $S\subset X$. The configurations $C_p$, $C'_p$
agree,
$u_i=u'_i$, except that $u_{p-1}=u'_{p-1}+S$. However, since $S$ is orthogonal
to all
the basic classes of $X$, for all $i$, all intersections of $u_i$ and $u'_i$
with all
basic classes of $Y=X\#(p-1)\CPb$ agree. Thus, according to
Theorem~\ref{BASIC}, the
coefficients $b_J$ are the same coefficients that arise in the formula
\[ {\D}_{X_p}=\exp(Q_{X_p}/2){\K}_X\sum_Jc_Je^{|J|S/p}. \] This means that
the sum of the coefficients of the expression for ${\D}_{X_p}$ in
Proposition~\ref{formallog} is $\sum_Jb_J=p$. \end{pf}

We next invoke Proposition ~\ref{ponq} to see that if $p$ is any positive odd
integer,
then a multiplicity $2p$ log transform can be obtained as the result of either
a
multiplicity $p$ log transform on a nodal fiber of multiplicity 2, or by a
multiplicity 2 log transform on a nodal fiber of multiplicity $p$. Thus
\begin{eqnarray*}
{\D}_{E(n;2p)}&=&\exp(Q/2)(e^{f_2}+e^{-f_2})(b_{p,0}+\sum_{i=1}^{(p-1)/2}
b_{p,2i}(e^{2if_2/p}+e^{-2if_2/p}))\\
&=&\exp(Q/2)(b_{p,0}+\sum_{i=1}^{(p-1)/2}b_{p,2i}(e^{2if_p}+e^{-2if_p}))
(e^{f_p/2}+e^{-f_p/2})
       \end{eqnarray*} since we already know the formula for a log transform of
multiplicity 2. We compare coefficients using $f_2=pf_{2p}$ and
$f_p=2f_{2p}$.

Assume for the sake of definiteness that $p\equiv 1\pmod4$ and let $r=(p-1)/4$.
In the
top expansion, the coefficient of $e^{\pm pf_{2p}}$ is $b_{p,0}$ and
$b_{p,2j}$ is the coefficient of $e^{\pm (p+2j)f_{2p}}$ and $e^{\pm
(p-2j)f_{2p}}$.
In the second expansion, the coefficient of $e^{\pm f_{2p}}$ is $b_{p,0}$, and
$b_{p,2j}$ is the coefficient of $e^{\pm (4j-1)f_{2p}}$ and $e^{\pm
(4j+1)f_{2p}}$. To
simplify notation, let $(m)_1$ be the coefficient of $e^{mf_{2p}}$ in the top
expansion
and $(m)_2$ its coefficient in the bottom expansion. Then,
\begin{eqnarray*}
b_{p,0}=(p)_1&=&(p)_2=b_{p,2r}=(p-2)_2=(p-2)_1=b_{p,2}=(p+2)_1=(p+2)_2
\\&=&b_{p,2(r+1)}=
(p+4)_2=(p+4)_1=b_{p,4}=(p-4)_1=(p-4)_2\\&=&b_{p,2(r-1)}=(p-6)_2=(p-6)_1=b_{p,6}=\cdots
\end{eqnarray*} and we see inductively that when $p$ is odd, all the $b_{p,2i}$
are
equal. But by Lemma~\ref{sum}
\[ b_{p,0}+2\sum_{i=1}^{(p-1)/2}b_{p,2i}=p.\] It follows that each
$b_{p,2i}=1$,
$i=0,\dots,(p-1)/2$.

Similarly, if $p$ is even, let $q=p-1$. Expanding ${\D}_{E(n;pq)}$ we see that
all
$b_{p,2i-1}$, $i=1,\dots,p/2$ are equal; and so again each $b_{p,2i-1}=1$.

\begin{thm}\label{toplt} Let $X$ be a $4$-manifold of simple type and suppose
that $X$
contains a nodal fiber $S$ orthogonal to all its basic classes. Then
\[ {\D}_{X_p}=\exp(Q_{X_p}/2){\K}_X\cdot{\sinh(S)\over\sinh(S/p)}. \] \end{thm}
\begin{pf} If, e.g., $p$ is odd,
\begin{eqnarray*} {\D}_{X_p}&=&
\exp(Q_{X_p}/2){\K}_X\cdot(1+2\cosh(2S/p)+2\cosh(4S/p)+\cdots+2\cosh((p-1)S/p))\\&=&
\exp(Q_{X_p}/2){\K}_X\cdot{\sinh(S)\over \sinh(S/p)}. \end{eqnarray*}  \end{pf}

As a result we have the calculation of the Donaldson series for all simply
connected
elliptic surfaces with $p_g\ge1$.

\begin{thm}\label{ellformula} If $n\ge2$ and $p,q\ge1$ are relatively prime,
\[{\D}_{E(n;p,q)}=\exp(Q/2){\sinh^n(f)\over\sinh(f_p)\sinh(f_q)}.\ \ \
\qed\]\end{thm}
\noindent This formula was originally conjectured by Kronheimer and Mrowka
\cite{KM}.

As an example of Theorem~\ref{toplt} consider $E(n)$. It follows from
\cite{GM1} and
\cite{FScusp} that in $E(n)$ there are 3 pairs of disjoint nodal fibers such
that the
nodal fibers in each pair are homologous, but give three linearly independent
homology
classes. Form $E(n;p_1,q_1;p_2,q_2;p_3,q_3)$ by performing log transforms with
each
pair $\{ p_i,q_i\}$ relatively prime. The resulting manifold is simply
connected and,
\begin{prop}\label{noncomplex}\hspace{.1in}$\displaystyle
{\D}_{E(n;p_1,q_1;p_2,q_2;p_3,q_3)}=
\exp(Q/2){\sinh^{n+4}(f)\over\prod\limits_{i=1}^3\sinh(f_{p_i})\sinh(f_{q_i})}.\ \
\qed $\end{prop}
\noindent Applying Theorem~\ref{ellformula} and
Proposition~\ref{noncomplex} to the manifolds $E(n; p_1,q_1;p_2,q_2; p_3,q_3)$,
we see
that they do not admit complex structures with either orientation
(cf.\cite{GM1},\cite[Theorem 8.3]{FScusp}).  The manifolds
$E(2;p_1,q_1;p_2,q_2;p_3,q_3)$ are the Gompf-Mrowka fake K3-surfaces
\cite{GM1}.

\bigskip

\section{Tautly Embedded Configurations\label{taut}}

Consider a $4$-manifold $X$ of simple type containing the configuration $C_p$.
By
Theorem~\ref{FSadj} for each 2-sphere $u_i$ in $C_p$ and each basic class
$\k$ of $X$, we have
\begin{equation}\label{tautly} -2\ge u_i^2+|u_i\cdot\k|
\end{equation} except in the special case described in
Theorem~\ref{FSadjspecial} where
$0\ge u_i^2+|u_i\cdot\k|$. The only examples known where the special case
arises are in
blowups. This was the situation in the previous section where we studied log
transforms. In this section, we assume  that we are not in the special case. We
say
that a configuration is {\em tautly embedded} if \eqref{tautly} is satisfied
for each
$u_i$ of the configuration and each basic class $\k$ of $X$.  Thus, if $C_p$ is
tautly
embedded, then for every basic class $\k$, $u_i\cdot\k=0$ for  $i=1,\dots,p-2$
and
$|u_{p-1}\cdot\k|\le p$.

\begin{thm}\label{tautcalc} Suppose that $X$ is of simple type and contains the
tautly
embedded configuration $C_p$. If
\[ {\D}_X=\exp(Q_X/2)\sum a_se^{\k_s} \] then the rational blowdown $X_p$
satisfies
\[ {\D}_{{X_p}}=\exp(Q_{X_p}/2)\sum \bar{a}_se^{\bar{\k}_s} \] where
\[\bar{a}_s=\begin{cases} 2^{p-1}a_s, \ \ \ &|u_{p-1}\cdot\k_s|=p\\
                           0, &|u_{p-1}\cdot\k_s|<p \end{cases}\] Furthermore,
if
$|u_{p-1}\cdot\k_s|= p$, then $\bar{\k}_s^2={\k_s}^2+(p-1)$.
\end{thm}
\begin{pf} If $\k_s\cdot u_{p-1}\ne 0,\pm p$ then $\bar{a}_s=0$ by
Lemma~\ref{theyextend}. For $\k_s\cdot u_{p-1}=0$, $\k_s\ne0$, note that since
the $\k_s$ are
characteristic, $p$ must be even.  But then  $\bar{\k}_s$ cannot even be
characteristic in $X_p$, since $\bar{\k}_s^2=\k^2$ is not mod $4$ congruent to
$(3\text{sign}+2e)(X_p)$.
Thus, Theorem~\ref{BASIC} implies
that $\bar{a}_s=0$.

In case $\k_s\cdot u_{p-1}=\pm p$, we compare with the model for the order $p$
log
transforms of $E(2)$; $C''_p\subset Y=E(2)\#(p-1)\CPb$ which is blown down to
obtain
$Y_p=E(2;p)$. Again let
$\lam_0=\pm(e_1+\cdots+e_{p-1})$; so by Lemma~\ref{J}, $\pm\lam_0$ are the
unique basic
classes of $Y_p$ satisfying $\pm\bar{\lam}_0=\pm(p-1)f_p\in H_2(Y_p;{\Z})$. Now
\begin{eqnarray*} {\D}_{Y}&=&{1\over2^{p-1}}\exp(Q/2)\sum\exp(\pm
e_1\pm\cdots\pm
e_{p-1})=
     \exp(Q/2)\sum_J{1\over2^{p-1}}e^{\lam_J}\\
{\D}_{Y_p}&=&\exp(Q/2)\sum\begin{Sb}|\l|\le p-1\\\l\equiv p\pmod2\end{Sb}e^{\l
f_p}=
     \exp(Q/2)\sum{1\over2^{p-1}}b_Je^{\lam_J}
\end{eqnarray*} Since $\pm\lam_0$ are the unique $\lam_J$ with
$\bar{\lam}_0=\pm
(p-1)f_p$, the corresponding coefficient is
$b_0=2^{p-1}$. We may now apply Theorem~\ref{BASIC} to obtain our result since
$\k_s\cdot u_i=\lam_0\cdot u''_i$ for each $i$.

In order to compute $\bar{\k}_s^2$, we find $x_i\in{\Q}$, $i=1,\dots,p-1$, such
that
\[ \k_s+\z=\k_s+\sum_{i=1}^{p-1}x_iu_i \in H_2(X^*;{\Q}) \] as in the proof of
Proposition~\ref{J}. We can solve for the $x_i$ using the model
$C''_p\subset Y''$, and $\pmb{\e}_J=\pm(1,\dots,1)$ in the proof of
Proposition~\ref{J}. Referring there, we get
\[ {\bold{x}}=P^{-1}A\pmb{\e}_J=-(A^t)^{-1}\pmb{\e}_J=\pm {1\over
p}(1,2,\dots,p-1). \]
So $\z=\pm\sum{i\over p}u_i$, and $\z^2={\bold{x}}\cdot P{\bold{x}}=1-p$.
Hence
\[ \bar{\k}_s^2 = (\k_s+\z)^2 = \k_s^2 + 2\k_s\cdot\z +\z^2= \k_s^2 +(p-1). \]
\end{pf}

Now consider the elliptic surface $E(1)$. It can be constructed by blowing up
${\bold{CP}}^2$ at the nine intersection points of a generic pencil of cubic
curves.  The
fiber class of $E(1)$ is $f=3h-e_1-\cdots e_9$ where $3h$ is the class of the
cubic in
$H_2({\bold{CP}}^2;{\Z})$. The nine exceptional curves are disjoint sections of
the
elliptic fibration. The elliptic surface $E(n)$ can be obtained as the fiber
sum of $n$
copies of $E(1)$, and these sums can be made so that the sections glue together
to give
nine disjoint sections of $E(n)$, each of square $-n$. In particular, consider
$E(4)$ with 9 disjoint sections of square $-4$. The basic classes of $E(4)$ are
$0$ and
$2f$; so we see that each of the 9 sections gives us a tautly embedded
configuration
$C_2$. Let $W_n$ be the rational blowdown of $n$ of these sections, $1\le n\le
9$.  For
$n\le 8$, $W_n$ is simply connected. Gompf has shown that all these manifolds
admit
symplectic structures, and it is not hard to see that $W_2$ is the 2-fold
branched
cover of ${\bold{CP}}^2$ branched over the octic curve \cite[\S5.2]{Gompf}.

\begin{prop} \hspace{.1in}$\displaystyle
{\D}_{W_n}=2^{n-1}\exp(Q/2)\cosh(\k_n)$ \
where \  $\k_n^2=n$.\end{prop}
\begin{pf} We have
\[ {\D}_{E(4)}=\exp(Q/2)\sinh^2(f)=\exp(Q/2)(\frac12\cosh(2f)+\frac12). \] The
basic
classes $\pm 2f$ intersect each section twice; so Theorem~\ref{tautcalc}
implies that
each $X_n$ has only the basic classes, $\pm\k_n$, and that each blowdown
multiplies its
coefficient by 2 and increases its square by $1$. (We start with coefficient
$\frac12$ and square $0$.)
\end{pf}

To further illustrate the utility of Theorem~\ref{tautcalc} we compute the
Donaldson
invariants of a family of Horikawa surfaces $\{ H(n)\}$ with
$c_1(H(n))^2=2n-6$. To obtain $H(n)$, start with the simply connected ruled
surface
${\bold{F}}_{n-3}$ whose negative section $s_-$ has square $-(n-3)$. We have
seen in
the proof of Lemma~\ref{ratball} that the classes $s_++f$ and $s_-$ form a
configuration in ${\bold{F}}_{n-3}$ whose regular neighborhood $D_{n-2}$ has
complement
the rational ball $B_{n-2}$. The Horikawa surface $H(n)$ is defined to be the
$2$-fold branched cover of ${\bold{F}}_{n-3}$ branched over a smoothing of
$4(s_++f)+2s_-$. (Equivalently this is a smooth surface representing (6,$n+1$)
in
$S^2\x S^2$.)

\begin{lem} For $n\ge 4$, the elliptic surface $E(n)$ contains a pair of
disjoint
configurations $C_{n-2}$ in which the spheres $u_{n-1}$  are sections of $E(n)$
and for
$1\le j\le n-2$, $u_j\cdot f=0$. Furthermore, the rational blowdown of this
pair of
configurations is the Horikawa surface $H(n)$.
\end{lem}
\begin{pf} It follows from our description of $H(n)$ that there is a
decomposition
\[ H(n)= B_{n-2}\cup \tilde{D}_{n-2}\cup B_{n-2} \] where $\tilde{D}_{n-2}$ is
the
branched cover of $D_{n-2}$.  Rationally blow up each
$B_{n-2}$; this is then the $2-$fold branched cover of ${\bold{F}}_{n-3}$ with
$B_{n-2}$ blown up. The result is the complex surface  $ C_{n-2}\cup
\tilde{D}_{n-2}\cup C_{n-2}$ which, by computing characteristic numbers, is
just $E(n)$.
\end{pf}
\noindent The first case, $n=4$, gives the example $H(4)=W_2$ above. The
Horikawa
surfaces $H(n)$ lie on the Noether line $5c_1^2-c_2+36=0$, and of course the
elliptic
surfaces $E(n)$ lie on the line $c_1^2=0$ in the plane of coordinates
$(c_1^2,c_2)$.
Let $Y(n)$ be the simply connected $4$-manifold obtained from $E(n)$ by blowing
down
just one of the configurations $C_{n-2}$. Then $c_1(Y(n))^2=n-3$ and
$c_2(Y(n))=11n+3$; so $Y(n)$ lies on the bisecting line $11c_1^2-c_2+36=0$.
The
calculation of Donaldson invariants of $Y(n)$ and $H(n)$ follows directly from
Theorem~\ref{tautcalc}.

\begin{prop} The Donaldson invariants of $Y(n)$ and $H(n)$ are:
\begin{eqnarray*} {\D}_{Y(n)}&=&\begin{cases}\exp(Q/2)\sinh(\lam_n),\ \ n\
\text{odd}\\
  \exp(Q/2)\cosh(\lam_n),\ \ n\ \text{even}\end{cases}\\
{\D}_{H(n)}&=&\begin{cases}2^{n-3}\exp(Q/2)\sinh(\k_n),\ \ n\ \text{odd}\\
  2^{n-3}\exp(Q/2)\cosh(\k_n),\ \ n\ \text{even}\end{cases}
\end{eqnarray*} where $\lam_n^2=n-3$ and $\k_n^2=2n-6$.\ \ \qed\end{prop}

\begin{cor} The simply connected $4$-manifolds $Y(n)$ are not homotopy
equivalent to
any complex surface.\end{cor}
\begin{pf} If  $Y(n)$ were homeomorphic to a complex surface, this computation
shows
that it  would have to be minimal, since the formula for ${\D}_{Y(n)}$ does not
contain
a factor $\cosh(e)$ where $e^2=-1$. Certainly the surface in question could not
be
elliptic since $c_1(Y(n))^2\ne 0$. Bt neither could the surface be of general
type
since $Y(n)$ violates the  Noether inequality. Thus $Y(n)$ is  not homeomorphic
to any
complex surface.\end{pf}

D. Gomprecht \cite{Gomprecht} has computed the value of the Donaldson invariant
$D_X(F^k)$ for any Horikawa surface $X$ and $k$ large, where $F$ is the
branched cover
of the fiber $f$ of $F_{n-3}$.

\bigskip

\section{Seiberg-Witten Invariants of Rational Blowdowns\label{SW}}

Suppose we are given a spin$^{\text{c}}$ structure on an oriented closed
Riemannian $4$-manifold $X$. Let $W^+$  and $W^-$ be the associated
spin$^{\text{c}}$ bundles with
$L=\det W^+=\det W^-$ the associated determinant line bundle. Since $c_1(L)\in
H^2(X;{\Z})$ is a characteristic cohomology class, i.e. has mod 2 reduction
equal to
$w_2(X)\in H^2(X;{\Z}_2)$, we refer to $L$ as a characteristic line bundle.
We will confuse a characteristic line bundle $L$ with its first Chern class
$L \in H^2(X;{\Z})$. For simplicity we assume that
$H^2(X;{\Z})$ has no $2$-torsion so that the set $Spin^{\text{c}}(X)$  of
spin$^{\text{c}}$  structures on $X$ is precisely the set of characteristic
line bundles on $X$.

Clifford multiplication, $c$, maps $T^\ast X$ into the skew adjoint
endomorphisms of
$W^+\oplus W^-$ and is determined by the requirement that
$c(v)^2$ is multiplication by $-|v|^2$. Thus $c$ induces a map
$$c: T^\ast X \to {\text{Hom}}(W^+, W^-).$$
 The $2$-forms $\Lambda^2=\Lambda^+\oplus\Lambda^-$ of $X$ then act on
$W^+$ leading to a map $\rho:\Lambda^+\to {\text{su}}(W^+)$.   A connection
$A$ on $L$ together with the Levi-Civita connection on the tangent bundle of
$X$ induces a connection
$\nabla_A:\Gamma(W^+)\to \Gamma(T^\ast X\otimes W^+)$ on $W^+$.  This
connection, followed by Clifford multiplication, induces the Dirac operator
$D_A:\Gamma(W^+)\to\Gamma(W^-)$. (Thus $D_A$ depends both on the connection
$A$ and the Riemannian metric on $X$.)  Given a pair $(A,\psi) \in
{\cal{A}}_X(L)\times
\G(W^+)$, i.e. $A$ a connection in $L$  and  $\psi$ a section of  $W^+$,  the
monopole equations of Seiberg and  Witten \cite{Witten} are
\begin{eqnarray}\label{monopole} D_A\psi&=&0\\\rho(F_A^+)\notag &=
&(\psi\otimes\psi^\ast)_o \end{eqnarray}
 where $(\psi\otimes\psi^\ast)_o$ is the trace-free part of the endomorphism
$\psi\otimes\psi^\ast$.

The gauge group
$\text{Aut}(L)=\text{Map}(X,S^1)$ acts on the space of solutions,  and its
orbit space is the moduli space $M_X(L)$ whose formal dimension is
\begin{equation} \dim M_X(L) = \frac14(c_1(L)^2-(3\,\text{sign}(X)+2\,e(X)).
 \label{dims} \end{equation}  If this formal dimension is nonnegative and if
$b^+>0$, then for a generic metric on $X$ the moduli space
$M_X(L)$ contains no reducible solutions (solutions of the form $(A,0)$
where $A$ is an \asd connection on $L$), and for a generic perturbation of the
second equation of \eqref{monopole} by the addition of a self-dual 2-form of
$X$, the moduli space $M_X(L)$ is a compact manifold of the given dimension
(\cite{Witten}).

The {\em Seiberg-Witten invariant} for $X$ is the function
$SW_X:Spin^{\text{c}}(X)\to {\Z}$ defined as follows. Let $L$ be a
characteristic line bundle.
  If $\dim M_X(L)<0$ or is odd, then
  $SW_X(L)$  is defined to be $0$. If $\dim M_X(L)= 0$, the moduli space
$M_X(L)$ consists of a finite collection of points and
$SW_X(L)$  is defined to be the number of these points counted with signs.
These signs are determined by an orientation on $M_X(L)$,  which in turn is
determined by an orientation on the determinant line
 $\det(H^0(X;{\bold{R}}))\otimes\det(H^1(X;{\bold{R}}))\otimes
\det(H^2_+(X;{\bold{R}}))$. If $\dim M_X(L)>0$ then we consider the basepoint
map
$$
\tilde{M}_X(L)=\{\text{solutions}\, (A,\psi)\}/\text{Aut}^0(L)\to M_X(L)
$$ where $\text{Aut}^0(L)$ consists of gauge transformations which are the
identity on the fiber of $L$ over a fixed basepoint in $X$. If there are no
reducible solutions, the basepoint map is an $S^1$ fibration, and we denote
its euler class by $\b\in H^2(M_X(L);\Z)$. The moduli space $M_X(L)$
represents an integral cycle in the configuration space
$B_X(L) =({\cal{A}}_X(L)\x \G(W^+))/\text{Aut}(L)$, and if $\dim M_X(L)=2n$,
the Seiberg-Witten invariant is defined to be the integer
$$
SW_X(L)=\la\b^n,[M_X(L)]\ra.
$$
A fundamental result is that if $b^+(X)\ge2$, the map
$$ SW_X: Spin^{\text{c}}(X) \to {\Z}
$$ is a diffeomorphism invariant ({\cite{Witten}); i.e.
 $SW_X(L)$ does not depend on the (generic) choice of Riemannian metric on
$X$ nor the choice of generic perturbation of the second equation of
\eqref{monopole}.

It is often convenient to observe that the space ${\cal{A}}_X(L)\x \G(W^+)$ is
contractible  and $\text{Aut}(L)\cong\text{Map}(X,S^1)$  acts freely on
${\cal{A}}_X(L)\x (\G(W^+)\setminus\{ 0 \})$.  Since $S^1$ is a $K({\Z},1)$, if
we
further assume that
$H^1(X;\R)=0$, then the quotient
$$ B^*_X(L)=\left({\cal{A}}_X(L)\x(\G(W^+)\setminus\{ 0 \})\right)/S^1
$$ of this action is homotopy equivalent to ${\bold{CP}}^{\infty}$. So if
there are no reducible solutions, we may view $M_X(L)\subset
{\bold{CP}}^{\infty}$. Under these identifications, the class
$\b$ becomes the standard generator of $H^2({\bold{CP}}^{\infty};{\Z})$.

Call a characteristic line bundle with nontrivial Seiberg-Witen invariant a
{\it Seiberg-Witten class}. The assumption in Seiberg-Witten theory which is
analogous
to the  assumption of simple type in Donaldson theory is
\begin{enumerate}\item[(*)] For each Seiberg-Witten class $L$, $\dim M_L(X)=0$.
\end{enumerate}
If this condition is satisfied, $X$ is said to have {\it Seiberg-Witten simple
type}.

\begin{lem}\label{char} Let $C_p\subset X$ and let $X_p$ be its rational
blowdown.
Assume that $X_p$ is simply connected.
\begin{enumerate}
\item[(a)] A line bundle $L^*$ on $X^*$ extends over $X_p$ if and only if
$c_1(L^*|_{L(p^2,1-p)})\in p{\Z}_{p^2}$.
\item[(b)] If $\bL$ is a characteristic line bundle on $X_p$, then there is a
characteristic line bundle $L$ on $X$ such that $L|_{X^*}=\bL|_{X^*}$.
\item[(c)] If $L$ is a characteristic line bundle on $X$, then an extension
$\bL$ of
$L|_{X^*}$ is characteristic on $X^*$ if and only if $\bL|_{B_p}$ is
characteristic.
\end{enumerate}\end{lem}
\begin{pf} (a) is obvious.
For any simply connected manifold $Y=V\cup W$ where $\bd V=\bd W$ is a rational
homology sphere, a class $c\in H^2(Y;\Z)$ will be characteristic provided
$\la c,\a\ra\equiv\a\cdot\a$ (mod 2) for all
$\a\in H^2(Y;\Z)$. Thus we need not worry about torsion classes; so a class is
characteristic if and only if its restrictions to $V$ and $W$ are both
characteristic. Applying this observation to $X_p=X^*\cup B_p$ proves (c).

To prove (b), let $\bL$ be a characteristic line bundle on $X_p$ and
let $L^*=\bL|_{X^*}$.  By (a), $\delta c_1(L^*)=mp$ for some
integer $m$. Suppose that $p$ is odd, then since $mp=(p+m)p\in\Z_{p^2}$, we may
assume
that $m$ is even. Let $L'$ be the line bundle on $C_p$ such that the Poincar\'e
dual
of $c_1(L')$ is $(m+1)\g_{p-1}+(m-p+1)\g_1$. Then $L'$ is characteristic on
$C_p$
and $\delta c_1(L')=\delta c_1(L^*)$. It follows that $L^*$
extends to a characteristic  line bundle on $X$ by our observation above. If
$p$ is
even, we may take $c_1(L')$ to be the Poincar\'e dual of $mp\g_1$ and get the
extension of $L^*$ to a global line bundle $L$ on $X$ whose restriction to both
$X^*$
and $C_p$ is characteristic.
\end{pf}

If $\bL$ is a line bundle on $X_p$ and $L$ is a line bundle on $X$ satisfying
$L|_{X^*}=\bL|_{X^*}$, we say that $L$ is a {\em lift} of $\bL$.

\begin{thm}\label{swgen} Let $C_p\subset X$ and let $X_p$ be its rational
blowdown.  Let $\bL$ be
a characteristic line bundle on $X_p$ and let $L$ be any lift of $\bL$ which is
characteristic on
$X$. Suppose that $\dim M_X(L)\equiv\dim M_{X_p}(\bL)\pmod 2$. Then
\[ SW_{X_p}(\bL)=SW_X(L).\] \end{thm}
\begin{pf} Since the rational ball $B_p$ embeds in the ruled surface
${\bold{F}}_{p-1}$ (see Lemma~\ref{ratball}), it admits a metric of positive
scalar curvature.
The gluing theory for solutions of the Seiberg-Witten equations follows the
same pattern as for
solutions of the \ASD equations. Thus we study the solutions on $X_p$ for $\bL$
by stretching the
neck between $X^*$ and $B_p$. We may assume that there are positive scalar
curvature metrics on
both the neck $L(p^2,1-p)\x\R^+$ and on $B_p$. This means that the only
solution to the
Seiberg-Witten equations on $B_p$ with a cylindrical end is the reducible
solution $(A',0)$,
where $A'$ is an \asd connection on $L'=\bL|_{B_p}$. Possible global solutions
are constructed
from asymptotically reducible  solutions on $X^*$ glued to $(A',0)$. The formal
dimension of
$M_{B_p}(L')$ is odd and negative, and there is one gluing parameter (since the
asymptotic value
is reducible); so
\[ \dim M_{X^*}(L^*)+1+\dim M_{B_p}(L') = \dim M_{X_p}(\bL)=2d_{\bL}, \] say.
(If $\dim M_{X_p}(\bL)$ is odd, there is nothing to prove.) Thus $\dim
M_{X^*}(L^*)=2d_{L^*}$
where $d_{L^*}\ge d_{\bL}$. This means that there is an obstruction to
perturbing a glued-up
$(A^*,\psi^*)\# (A',0)$ to a solution. As in Donaldson theory, there is an
obstruction bundle
$\xi$ over $M_{X^*}(L^*)$, and it is the  complex vector bundle of rank
$d_{L^*}-d_{\bL}$
associated to the basepoint fibration. The zero set of a generic section of
$\xi$ is homologous
to $M_{X_p}$ in $B_{X_p}(\bL)$. Thus
\[ SW_{X_p}(\bL)=\la\b^{d_{\bL}},[M_{X_p}(\bL)]\ra=\la\b^{d_{\bL}},
\b^{d_{L^*}-d_{\bL}}\cap [M_{X^*}(L^*)]\ra=\la\b^{d_{L^*}},[M_{X^*}(L^*)]\ra.
\]

Let $L$ be a characteristic line bundle on $X$ which is a lift of $\bL$, and
let
$\dim M_X(L)=2d_L$. The second construction of
Lemma~\ref{ratball} shows that $C_p$ has a metric of positive scalar curvature.
So the discussion
of the last paragraph applies to show that
\[ SW_X(L)=\la\b^{d_{L^*}},[M_{X^*}(L^*)]\ra, \]
completing the proof of the theorem.
\end{pf}

\begin{lem}\label{swdim} Suppose that $C_p\subset X$ with rational blowdown
$X_p$.
Let $L$ be a characteristic line bundle on $X$ such that $\la c_1(L), u_i\ra=0$
for
$i=1,\dots,p-2$ and $\la c_1(L), u_{p-1}\ra=mp$ for some $m\in\Z$. Let $\bL$ be
a
characteristic extension of $L|_{X^*}$ to all of $X_p$. Then $m$ is odd, and
$\dim M_{\bL}(X_p)=\dim M_L(X)+{m^2-1\over4}(p-1)$.\end{lem}
\begin{pf} The proof of Theorem~\ref{tautcalc} shows that $m$ must be odd if
$\bL$ is to
be characteristic and that $c_1(\bL)^2=c_1(L)^2+m^2(p-1)$.  Since
$3\,\text{sign}(X_p)+2\,e(X_p)=3\,\text{sign}(X)+2\,e(X)+(p-1)$, the lemma
follows.
\end{pf}
\noindent Note that this shows that, unless
$m=\pm1$, the dimensions of the moduli spaces will increase.

We shall consider the two situations analogous to those studied in the previous
sections:
\smallskip
\begin{itemize}
\item[(i)] $C_p$ is embedded in $X=Y\#(p-1)\CPb$ so that $X_p$ is the result of
an
order $p$ log transform performed on a nodal fiber of $X$.\\
\vspace{-.1in}\item[(ii)] $C_p$ is tautly embedded in $X$ with respect to $L$,
i.e.
$\la c_1(L),u_i\ra=0$ for $i=1,\dots,p-2$, and $\la c_1(L),u_{p-1}\ra\le p$.
\end{itemize}

The next theorem follows directly from Theorem~\ref{swgen} and
Lemma~\ref{swdim}.
\begin{thm}\label{tautagain} Suppose that $X$ has Seiberg-Witten simple type
and that
$C_p\subset X$ with $X_p$ its rational blowdown. Assume that $X_p$ is simply
connected and that $\bL$ is a characteristic line bundle on $X_p$. Suppose
further that
$L$ is a characteristic lift of $\bL$ and that $C_p$ is tautly embedded with
respect to $L$. Then
\[ SW_{X_p}({\bL})=SW_X(L)\]
and $c_1(\bL)^2=c_1(L)^2+(p-1)$. \ \  \qed
\end{thm}

Say that the configuration $C_p$ is {\it SW-tautly embedded} in $X$ if it is
tautly
embedded with respect to each Seiberg-Witten class.

\begin{cor} Suppose that $X$ has Seiberg-Witten simple type and contains the
SW-tautly
embedded configuration $C_p$.  Assume that the rational blowdown $X_p$ is
simply
connected. Then the Seiberg-Witten classes of $X_p$ are the
characteristic line bundles $\bL$ which have a lift to a Seiberg-Witten class
$L$ of
$X$, and $SW_{X_p}({\bL})=SW_X(L)$. \ \ \qed\end{cor}

In a fashion similar to the proof of Theorem~\ref{swgen}, one can prove a
blowup formula
for Seiberg-Witten invariants. The characteristic line bundles of $X\#\CPb$ are
those of
the form $L\otimes E^{2k+1}$ where
$L$ is characteristic on $X$ and $c_1(E)=e$, and
$\dim M_{L\otimes E^{2k+1}}(X\#\CPb)=\dim M_L(X)-k(k+1)$. It is shown in
\cite{Turkey}
that $SW_{X\#\CPb}(L\otimes E^{2k+1})=SW_X(L)$ provided $\dim
M_L(X)-k(k+1)\ge0$. It
follows that if $X$ satisfies the Seiberg-Witten simple-type condition (*),
then so
does $X\#\CPb$.

Suppose that $X$ contains the nodal fiber $S$, and $X_p$ is the result of
performing
an order $p$ log transform on $S$. The characteristic line bundles on $X_p$ are
obtained from characteristic bundles
$L\otimes E_1^{2k_1+1}\otimes\cdots\otimes E_{p-1}^{2k_{p-1}+1}$ on
$Y=X\#(p-1)\CPb$ by
restricting to $Y^*=Y\setminus C_p$ and then extending over $B_p$. If we assume
that
$\la c_1(L),S\ra=0$, then for each $L\otimes E_1^{\pm1}\otimes\cdots\otimes
E_{p-1}^{\pm1}=L({\pmb\e}_J)$  with $c_1(L({\pmb\e}_J))=
c_1(L)+\sum_J\e_{J,i}e_i$, it
follows from Proposition~\ref{J} that the unique extension $\bL_J$ over $X_p$
has
$c_1(\bL_J)=c_1(L)+|J|\s_p$, where $\s_p$ is the Poincar\'e dual of $S/p$.
(Note
that when $p$ is even, $|J|$ must be odd; so the extension $\bL_J$ is
characteristic.)  Hence
\[ \dim M_{\bL_J}(X_p) =\dim M_L(X),\]
and Theorem~\ref{swgen} implies:

\begin{thm}\label{SWlog} Suppose that $X$ has Seiberg-Witten simple type and
contains the
nodal fiber $S$. Let $L$ be a characteristic line bundle on $X$ with $\la
c_1(L),S\ra=0$.
Let $X_p$  be the result of performing an order $p$ log transform on $S$. For
each $J\in
\{\pm1\}^{p-1}$, we have $SW_{X_p}({\bL_J})=SW_X(L)$.  Suppose furthermore that
$\la
c_1(L),S\ra=0$ for each characteristic $L$ on $X$ with $SW_X(L)\ne0$.  Then
$X_p$ also
has Seiberg-Witten simple type and each line bundle $\Lambda$ on $X_p$ with
$SW_{X_p}(\Lambda)\ne 0$ is of the form $\Lambda=\bL_J$.\ \ \qed
\end{thm}

By a the {\em nodal configuration} we shall mean a configuration
$C_p\subset X\#(p-1)\CPb$ as above, obtained from a nodal fiber $S$ satisfying
the
condition $\la c_1(L),S\ra=0$ for each characteristic $L$ on $X$ with
$SW_X(L)\ne0$.

Witten \cite{Witten} has conjectured that (for manifolds with
$b^+>1$) the Seiberg-Witten simple type condition is equivalent to the
simple type condition of Kronheimer and Mrowka for Donaldson theory. Further,
under this
hypothess of simple type, Witten gives a precise conjecture for relating the
Seiberg-Witten invariants and the Donaldson series, namely:

\begin{conj}[Witten]
The set of basic classes in the two theories are the same, and
\[ {\D}_X= 2^{3\text{sign}+2e-({b^+-3\over 2})}\exp(Q/2)\sum
SW_X(\k_s)e^{\k_s}.\]
\end{conj}

\begin{thm} Witten's conjecture is true for simply connected elliptic
surfaces.\end{thm}
\begin{pf} Witten has given a recipe for calculating $SW_X$ for all Kahler
manifolds
$X$. So one could prove this theorem simply by comparing the answer obtained
with that
of Theorem \ref{ellformula}. Alternatively, note that Witten's recipe gives the
result
that the nonzero Seiberg-Witten invariants of $E(n)$ are:
\begin{equation}\label{SWEn}
 SW_{E(n)}((n-2-2r)f)= (-1)^r\binom{n-2}{r},\ \ r=0,\dots,n-2
\end{equation}
 (where $f$ is the fiber class). Suppose we define
\[{\bold{W}}_X=2^{3\text{sign}+2e-({b^+-3\over 2})}\sum SW_X(\k_s)e^{\k_s}, \ \
{\bold{SW}}_X=\exp(Q_X/2){\bold{W}}_X\]
Then \eqref{SWEn} shows that ${\D}_{E(n)}={\bold{SW}}_{E(n)}$. Suppose
that $X_p$ is the result of an order $p$ log transform on a nodal fiber which
is
orthogonal to all classes in $H_2(X)$ with nontrivial Seiberg-Witten
invariants. Then
Theorem \ref{SWlog} implies that
${\bold{W}}_{X_p}={\bold{W}}_X\cdot(\sinh(f_p)/\sinh(f))$.
It follows that ${\bold{SW}}_{E(n;p,q)}={\D}_{E(n;p,q)}.$\end{pf}

Furthermore, we have

\begin{thm} If $X$ satisfies the Witten conjecture, then so do all blowups and
blowdowns and any rational blowdown $X_p$ of a nodal or taut configuration.
\end{thm}

\newpage
\centerline{\large{\sc{Figures}}}

\vspace*{3.15in}\hspace{1in} {\sc Figure 1.} \hspace{2in} {\sc Figure 2.}

\hspace*{1.2in} \hspace{2.5in}

\vspace*{3.15in}\hspace{1in}{\sc Figure 3.} \hspace{2in} {\sc Figure 4.}

\hspace*{1in} \hspace{2.35in}
\newpage

\vspace*{3.15in}\hspace{1in} {\sc Figure 5.} \hspace{2in} {\sc Figure 6.}

\hspace*{.9in}

\vspace*{3.15in}\hspace{1in} {\sc Figure 7.} \hspace{2in} {\sc Figure 8.}

\newpage

\vspace*{3.15in}\hspace{1in} {\sc Figure 9.} \hspace{2in} {\sc Figure 10.}

\hspace*{.9in}

\vspace*{3.15in}\hspace{1in} {\sc Figure 11.} \hspace{2in} {\sc Figure 12.}

\newpage

\vspace*{3.15in}\hspace{1in} {\sc Figure 13.} \hspace{2in} {\sc Figure 14.}

\end{document}